%% file: main.tex
\documentclass{article}
\usepackage{xcolor}
\usepackage{keyval}
\usepackage{constants}
\newconstantfamily{epsilons}{
symbol=\varepsilon
}
\input{1forms_preamble.tex}
\usepackage{amsmath}
\usepackage{todonotes}
\author{Arham Deep}
\date{\today}
\title{Explicit elliptic estimates for nowhere vanishing harmonic $1$-forms}
\definecolor{amethyst}{rgb}{0.6, 0.4, 0.8}
\usepackage{natbib}
\begin{document}
\maketitle
\begin{abstract}
We compute an explicit constant for an injectivity estimate on $T^3=S^1 \times S^1 \times S^1$ involving the Laplace Operator. First, we provide motivation for such explicit estimates. We perform the computation for $T^3$ endowed with the flat metric $\gflat$ before generalising to perturbed metrics. Finally, we apply these results to show existence of a nowhere vanishing harmonic 1-form on the 3-Torus endowed with a perturbed metric. 
\end{abstract}
\tableofcontents

\section{Introduction}
The Laplace operator on a compact manifold $M$ is one of the best understood elliptic differential operators, and much is known about it.
In particular, \emph{standard elliptic theory}, which was developed over the course of the last century, provides the following \emph{injectivity estimate}:

\begin{theorem}[Theorem H.27 in \cite{Besse1987}]
    There exists a constant $c>0$ depending on $M$ (but not on $f$) such that the following is true:
    for all $f \in W^{2,p}(M)$ which are $L^2$ -orthogonal to $\Ker \Delta$ we have that:
    \begin{align}
    \label{equation:standard-elliptic-estimate}
        \|{f}_{W^{2,p}}
        \leq
        c \|{\Delta f}_{L^p}.
    \end{align}
\end{theorem}

This theorem is typically shown by contradiction, and the constant $c$ depends on the geometry of the manifold in a highly complicated way.
In many problems in analysis, existence of the constant $c$ alone is not enough, but it is necessary to give an \emph{upper bound} for it.

Often, the manifold $M$ is given through an explicit construction depending on some parameter, and it suffices to understand how $c$ (or an upper bound for it) asymptotically depends on this parameter for extreme values.
One example is the construction of anti-self-dual instantons in \cite{Taubes1982}, pioneering this technique in geometry.
In this application, the parameter controlled the size of a connecting piece along which two manifolds are glued together.
Another example is the first counter-example to Payne's conjecture in \cite{Hoffmann1997}.
Here, the parameter controlled the number of holes of a region in $\R^2$.

However, even more challenging are cases in which the manifold $M$ does \emph{not} depend on a parameter, and it does not suffice to understand some asymptotic behaviour.
In these cases, one needs \emph{explicit elliptic estimates}.
Such estimates have been studied in their on right:
in \cite{Plum1992}, explicit estimates for a class of second-order operators were proved, depending on the shape of a region in $\R^2$.
In \cite{Guneysu2018,Guneysu2018b} estimates on manifolds, depending explicitly on geometric data such as the curvature, were proved.
Recently, the interest in such explicit estimates has been renewed, because they are required to carry out \emph{numerically verified proofs}.
\cite[Section 4.1]{Nakao2019} has examples of such proofs making use of explicit elliptic estimates.

Almost all explicit estimates appearing in the literature are for domains in $\R^n$ and for a weaker norm instead of the $W^{2,p}$-norm on the left side of \cref{equation:standard-elliptic-estimate}.
We are not aware of a work which gives an explicit constant for \cref{equation:standard-elliptic-estimate} on a closed manifold.
In this paper we prove such estimates:
\begin{itemize}
\item For $T^3$ endowed with the flat metric 
\item For $T^3$ endowed with a metric "close" to the flat metric.
\end{itemize}
The injectivity estimate for $T^3$ endowed with the flat metric follows readily from the local theory and we obtain the following result:

\begin{proposition}
\label{proposition:injectivity-estimate-flat-metric}
    For all $f \in W^{2,4}(T^3)$ which are $L^2$-orthogonal to $\Ker \Delta$ we have that:
    \begin{align}
        \|{f}_{W^{2,4}}
        \leq
        \Cl{const:injectivity-main-flat} \|{\Delta f}_{L^4},
    \end{align}
    where $\Cr{const:injectivity-main-flat}$ is defined in \cref{equation:flat-injectivity-constant}.
\end{proposition}

If one perturbs the metric, then the same estimate holds, only with a slightly perturbed constant.
That is made precise in the following result:

\begin{proposition}
    \label{theorem:T3-injectivity-estimate}
    Let $g_{\text{flat}}$ be the flat metric on $T^3$ and $g$ another metric.
    If $\|{g-g_{\text{flat}}}_{C^1} \leq \Cl{const:metric-max-difference}$, then for all $f \in L^p_2(T^3)$:
    \begin{align}
        \|{f}_{W^{2,4}}
        \leq
        \Cl{const:injectivity-main-nonflat} \|{\Delta f}_{L^4},
    \end{align}
    where $\Cr{const:injectivity-main-nonflat}$ is defined in \cref{equation:injectivity-nonflat-constant-def} and 
    $\Cr{const:metric-max-difference} $ must satisfy $\Cr{Absorption Laplacian Term}(\Cr{const:metric-max-difference}) < \frac{1}{\Cr{const:injectivity-main-flat}}$.
    Here, $\Cr{Absorption Laplacian Term}$ is defined in \cref{Delta Comparison Bound} and $\Cr{const:injectivity-main-flat}$ is defined in \cref{equation:flat-injectivity-constant}.
\end{proposition}

Our choice of $T^3$ is motivated by \cite[Section 7.5]{Joyce2021}:
there, a $3$-torus contained in a $6$-manifold is considered.
It inherits a metric from the ambient manifold whose definition is complicated, but it is expected to be close to (but not equal to) the flat metric on $T^3$.
If the metric was the flat metric, it would admit a harmonic $1$-form that is nowhere vanishing.
Since the metric on $T^3$ is conjectured to be close to be flat, it is still expected to admit a harmonic $1$-form that is nowhere vanishing.

We can use \cref{theorem:T3-injectivity-estimate} to give an effective bound for this.
Using it, we find the following \emph{sufficient criterion} for how close a metric on $T^3$ needs to be to the flat metric so that it admits a suitable $1$-form.
This makes \cite[Section 7.1(E)]{Joyce2021} more precise:

\begin{theorem}
    \label{theorem:application-nowhere-vanishing-1-form}
    Let $\gflat$ be the flat metric on $T^3$.
    If $\|{g-g_{\text{flat}}}_{g_{\text{flat}},C^1} \leq \Cl{const:metric-max-difference-for-one-forms}$, then there exists $\xi \in C^\infty(M)$ such that $\d x_1+ \d \, (\xi) \in \Omega^1(M)$ is a nowhere vanishing $1$-form that is harmonic with respect to $g$.
    Here, $\Cr{const:metric-max-difference-for-one-forms}$ must satisfy $\Cr{const:metric-max-difference-for-one-forms}<\Cr{const:metric-max-difference}$ and $\Cr{endofpaper}(\Cr{const:metric-max-difference-for-one-forms})>0$, where $\Cr{endofpaper}$ is defined in \cref{End of Proof}.
\end{theorem}

The article is structured as follows:
in \cref{subsection:Background} we list some basic definitions of geometric analysis to fix notations.
In \cref{section:explicit-elliptic-estimates} we prove the aforementioned injectivity estimates on $T^3$:
in \cref{subsection:local-estimates} we give explicit constants for some standard local estimates from the literature;
from this \cref{proposition:injectivity-estimate-flat-metric} quickly follows and the remainder of the proof is presented in \cref{subsection:estimates-on-t3-with-flat-metric};
then \cref{theorem:T3-injectivity-estimate} follows from comparing nearby metrics, which is done in \cref{subsection:estimates-on-t3-with-non-flat-metric}.
Last, in \cref{section:application-nowhere-vanishing-1-forms} we prove our application to nowhere vanishing harmonic $1$-forms.
\\
\textbf{Acknowledgements}:\\
The author would like to give special thanks to \textbf{Dr Daniel Platt} for supervising this project over the past year and a half, for taking the time to host regular meetings in order to discuss this project.\\
The author would also like to thank \textbf{Zihan Zhang} for keeping us company during the meetings and his great enthusiasm shown toward the project.
The author would like to thank the following people for many fruitful conversations about mathematics:
\begin{itemize}
\item Yiheng (Jackie) Wu
\item Rohan Money Shenoy
\item Ariff Jazlan Johan
\end{itemize}

\section{Explicit elliptic estimates}
\label{section:explicit-elliptic-estimates}
\subsection{Background}
\label{subsection:Background}
The standard definition of a Sobolev Space over a general manifold from \cite[2.1]{HebeyRobert2008} is the following:
\\
For Riemannian manifold $(M,g)$, the Sobolev Space $W^{k,p}(M,g)$ is the completion of $C^\infty (M)$ for the norm:
\begin{align*}
\|{u}_{W^{k,p}(M,g)} = \sum_{i=0}^k \|{\nabla^i u}_{g, L^p(M)},
\end{align*}
where $\|{\nabla^i u}_{g, L^p(M)}$ is the $L^p$-norm of the function $|\nabla^i u|$ with respect to g.
\\Here, we use the pointwise norm of a $(k,l)$ tensor $T$ as the norm defined by the following inner product
\begin{align}
\label{Induced inner product}
\langle F, G \rangle
= 
g^{i_1r_1}...g^{i_kr_k}g_{j_1s_1}...g_{j_ls_l} F^{j_1...j_l}_{i_1...i_k} G^{s_1...s_l}_{r_1...r_k}
\end{align}
For $M=\Omega \subset \R^N$, the pointwise norm of the $(0,1)$-tensor $\nabla u$ is equal to the vector norm of the gradient of $u$, and the pointwise norm of the $(0,2)$-tensor $\nabla^2 u$ is equal to the Frobenius norm of the Hessian of $u$.
\subsection{Local estimates}
\label{subsection:local-estimates}
In this section, we work in (open subsets of) Euclidean Space endowed with the flat metric. Any Sobolev norm in this subsection is assumed to be taken with respect to the flat metric. 

Let $\Omega \subset \R^n$ bounded and let $f \in L^p(\Omega)$ for some $1<p<\infty$.
Define the Newtonian potential of $f$ as
\[
\label{Newtonian Potential}
    w(x)
    :=
    \int_\Omega
    \Gamma(x-y)f(y) \, \d y,
\]
where $\Gamma$ is the fundamental solution of Laplace's equation in dimension $n$.

\begin{theorem}[Theorem 9.9 in \cite{Gilbarg1998}]
    \label{theorem:calderon-zygmund-potential-formulation}
    We have that $w \in L^p_2(\Omega)$ and $\Delta w=f$ almost everywhere, and
    \begin{align}
        \|{D_{ij} w}_{L^p}
        \leq
        C_{\text{Calderon-Zygmund}}(n,p)
        \|{f}_{L^p},
    \end{align}
    where $C_{\text{Calderon-Zygmund}}(n,p)=1$ if $p=2$ and for $p > 2$:
    \[
        C_{\text{Calderon-Zygmund}}(n,p)
        =
        C_{\text{Marcinkiewicz}}(1,2,p')
        \Cr{constant:mu-t-L1-norm} ^\alpha \Cr{constant:mu-t-L2-norm}^{1-\alpha},
    \]
    And for $1<p< 2 $:
    \[
        C_{\text{Calderon-Zygmund}}(n,p)
        =
        C_{\text{Marcinkiewicz}}(1,2,p)
        \Cr{constant:mu-t-L1-norm} ^\alpha \Cr{constant:mu-t-L2-norm}^{1-\alpha},
    \]
    where $\frac{1}{p}+\frac{1}{p'}=1$, the numbers $C_{\text{Marcinkiewicz}}(1,2,p')$ and $\alpha$ are defined in \cref{theorem:marcinkiewicz}.
\end{theorem}
\begin{proof}
We use the proof presented in \cite{Gilbarg1998}, which aims to use the Marcinkiewicz Interpolation Theorem (\cref{theorem:marcinkiewicz}) for the case $1<p<2$ and extend by duality. The Marcinkiewicz Interpolation Theorem gives us a sufficient condition to show that a linear map from $L^p(\Omega)$ into itself is bounded, namely two inequalities regarding the distribution functions of T (refer to \cref{theorem:marcinkiewicz})
\\
Define the linear operator: $T:L^2(\Omega) \rightarrow L^2(\Omega) $ as 
\begin{align}
\label{Definition of T}
Tf = D_{ij} w
\end{align}
Where $i,j$ are fixed and $w$ denotes the Newtonian Potential. 
\\
By the above discussion it is sufficient to show the following two inequalities:
\begin{align}
\label{Distribution Function Bounds}
\begin{split}
\mu_{Tf}(t) \leq \Cl{constant:mu-t-L2-norm} \left( \frac{\|{f}_{L^2}}{t} \right)^2,
\\
\mu_{Tf}(t) \leq \Cl{constant:mu-t-L1-norm} \frac{\|{f}_{L^1}}{t}.
\end{split}
\end{align}
Where $\mu_f(t):= |\{x \text{ such that } f(x)>t  \}|$.\\
$\Cr{constant:mu-t-L2-norm}=1$ is derived explicitly in \cite[Equation 9.29]{Gilbarg1998}.
\\We derive an explicit value for the constant $\Cr{constant:mu-t-L1-norm}$ and break down the derivation of this constant into steps. \\
\textbf{Step 1: Rewrite $f = b + g\in L^2(\Omega)$ as a sum of two functions for appropriate $b, g \in L^2(\Omega)$} \\
Our choice of functions $b,g \in L^2(\Omega)$ is motivated by the Cube Decomposition presented in \cite[Section 9.2]{Gilbarg1998}.
We first extend $f$ to vanish outside $\Omega$ and we fix a cube $K_0$ satisfying $\Omega \subset K_0$, chosen to be large enough such that for a fixed choice of $t>0$ we have that:
\begin{align*}
    \int_{K_0} |f| \leq t|K_0|.
\end{align*}
$|K_0|$ denotes the Lebesgue measure of the cube. Subdivide $K_0$ into $2^n$ congruent subcubes. We organise the subcubes as follows: let $K_2$ denote the collection of subcubes where each subcube $K \in K_2$ satisfies:
\begin{align}
\label{good}
    \int_{K}|f| \leq t |K|.
\end{align}
The remaining subcubes satisfy
\begin{align}
\label{bad}
    \int_{K}|f| > t |K|.
\end{align}
and we add them to the set $\Upsilon$. 
For each of the subcubes in $K_2$ we repeat the process and decompose it into $2^n$ subcubes. 
If they satisfy \cref{good}, then they are a part of a newly defined set $K_3$, if not we add them to the set $\Upsilon$. 
This process is continued inductively, and we denote by $\Upsilon_l$ an enumeration of all cubes in $\Upsilon$.
For each subcube $K \in \Upsilon$, we denote $\tilde{K}$ to be the subcube whose subdivision yields $K$. 
By construction, $|\tilde{K}|/ |K|= 2^n$.
So, for $K \in \Upsilon$ we find:
\begin{align*}
\int_{K} |f| \leq \int_{\tilde{K}} |f| \leq t |\tilde{K}| = 2^nt |K|.
\end{align*}
Combining this with our assumption that $K \in \Upsilon$ we obtain the following inequality.
\begin{align}
\label{nice property of bad functions}
    t < \frac{1}{|K|} \int_{K}|f| \leq 2^nt.
\end{align}
 By the Lebesgue Differentiation Theorem \cite[Chapter 3, Theorem 1.3]{SteinShakarchi2005} we have that, $|f| \leq t$ almost everywhere on $G := K_0 - \cup_l \Upsilon_l$.
 Now we split our function $f$ into a good and bad part, as done in \cite[p.232]{Gilbarg1998}.
 Define the good part, $g$, as follows:
\begin{align*}
    g(x)= \begin{cases}
        f(x) & \text{for $x \in G$}\\
        \frac{1}{|\Upsilon_l|}\int_{\Upsilon_l} f & \text{for $x \in \Upsilon_l$}.
    \end{cases}
\end{align*}
The bad part is defined as $b=f-g$, by the linearity of T we have that:
\begin{align}
\label{good+bad}
    \mu_{Tf}(t) \leq \mu_{Tg}(t/2)+ \mu_{Tb}(t/2).
\end{align} 
We list some properties of the "good" and "bad" function which we will use later \cite[p.226]{Gilbarg1998} :
\begin{align}
\label{good and bad properties}
\begin{split}
|g|\leq 2^n t \text{ a.e.},\\
b(x) =0 \text{  for x} \in G, \\
\int_{\Upsilon_l} b =0.
\end{split}
\end{align}
\textbf{Step 2: Show that there exists a constant $\Cl{constant:mu-tg-L1-norm}>0$ such that
\begin{align*}
\mu_{Tg} \left( \frac{t}{2} \right) \leq \frac{\Cr{constant:mu-tg-L1-norm}}{t} \|{f}_{L^1}. 
\end{align*}}\\
The explicit estimate for $\mu_{Tg}(t/2)$ is obtained in \cite[p.232]{Gilbarg1998}, this is:
\begin{align}
    \label{equation:mu-Tg-estimate}
    \mu_{Tg} \left( \frac{t}{2} \right) \leq \frac{2^{n+2}}{t} \|{f}_{L^1}.
\end{align}
We also provide the derivation here for completeness:
\begin{align*}
    \mu_{Tg} \left( \frac{t}{2} \right) &\leq \frac{4}{t^2}\int_{\Omega} g^2 \leq \frac{2^{n+2}}{t}\int_{\Omega} |g|.
\end{align*}

For the first inequality we use Markov's Inequality \cref{Markov Inequality} and for the second we use \cref{good and bad properties}

\textbf{Step 3: Show that there exists a constant $\Cl{constant:mu-tb-L1-norm}>0$ such that:}
\begin{align}
\label{very bad function estimate}
   \mu_{Tb}(\frac{t}{2}) \leq \frac{\Cr{constant:mu-tb-L1-norm}}{t} \|{f}_{L^1}.
\end{align}\\
As done in \cite[p.233]{Gilbarg1998} we define:
\begin{align}
    \label{equation:bl-definition}
    b_l = \begin{cases}
    b & \text{on $\Upsilon_l$},\\
    0 & \text{elsewhere}.
    \end{cases}
\end{align}
Hence we have that:
\begin{align*}
    Tb = \sum_{i=1}^{\infty} Tb_l.
\end{align*}

Now it suffices to show that there exists a constant such that for each $l$ we have the following:
\begin{align}
\label{step 3 Calderon-Zygmund}
    \left|
    x \in K_0 \text{   such that  } Tb_l(x) > \frac{t}{2}
    \right| 
    \leq \frac{\C}{t} \|{f}_{L^1}.
\end{align}
\begin{remark}
\label{idea of step 3 Calderon-Zygmund}
Define $\rho := \text{diam}(\Upsilon_l)$ and let $B_l$ be a ball of radius $\rho$ centred at $\bar{y}$, where $\bar{y}$ denotes the centre of $\Upsilon_l$. 
The rough idea is to obtain bounds of the form:
\begin{align}
\label{step 3 Calderon-Zygmund Hard}
    \left|
    x \in K_0 - B_l \text{   such that  } Tb_l(x) > \frac{t}{2}
    \right| 
    \leq \frac{\C}{t} \|{f}_{L^1}.
\end{align}
\begin{align}
    \left|
    x \in B_l \text{   such that  } Tb_l(x) > \frac{t}{2}
    \right| 
    \leq \frac{\C}{t} \|{f}_{L^1}.
\end{align}
Combining the above two bounds will give \cref{step 3 Calderon-Zygmund}. The choice to remove $B_l$ specifically is made for the following reason: we aim to use Markov's inequality to show \cref{step 3 Calderon-Zygmund Hard} and the corresponding integral can then be written in spherical coordinates. To obtain the latter bound, we bound the measure of the set $B_l$. 
\end{remark}
For fixed $l$ we approximate $b_l$ by a sequence of functions $\{b_{lm}\} \in C_0^{\infty}(\Upsilon_l)$ with compact support converging to $b_l$ in $L^2(\Omega)$.
For $x \notin \Upsilon_l$, we have the formula:
\begin{align*}
    Tb_{lm}(x) &= \int_{\Upsilon_l}{D_{ij}\Gamma(x-y) b_{lm}(y) dy}
               = \int_{\Upsilon_l} \{D_{ij}\Gamma(x-y) -D_{ij}\Gamma(x-\bar{y})\}b_{lm}(y) dy,
\end{align*}
Where we use the definition of $T$ \cref{Definition of T} in the first equality and \cref{good and bad properties} for the second.
We apply the mean value inequality: 
there exists $\tilde{y}$ between $y$ and $\bar{y}$ satisfying:
\begin{align*}
|Tb_{lm}(x)| \leq \int_{\Upsilon_l} |DD_{ij}\Gamma(x-\tilde{y})| \cdot |y-\bar{y}| |b_{lm}(y)dy.
\end{align*}
From the estimate of $I_6$ in the proof \cite[Lemma 4.4]{Gilbarg1998} we have
\begin{align}
    \label{lemma:2nd-Derivative-Newtonian}
    |DD_{ij}\Gamma(x-y)| \leq 
    \Cl{constant:-2nd-Derivative-Bound-Newtonian-Potential}
    |x-y|^{-n-1}.
\end{align}
Where $\Cr{constant:-2nd-Derivative-Bound-Newtonian-Potential}= \frac{n(n+5)}{\omega_n}$, consequently:
\begin{align*}
|Tb_{lm}(x)| \leq \Cr{constant:-2nd-Derivative-Bound-Newtonian-Potential}\delta[\operatorname{dist} (x, \Upsilon_l)]^{-n-1} \int_{\Upsilon_l} |b_{lm}(y)| dy.
\end{align*}

Recall $B_l$ is the ball centred at $\bar{y}$ of radius $\rho$, i.e. $B_l = B_{\delta}(\bar{y})$. 
Let $r(x):=\operatorname{dist} (x, \Upsilon_l)$. 
By integrating both sides and switching to spherical coordinates we arrive at:
\begin{align}
\label{equation:Tblm-estimate-on-K0-minus-Bl}
\begin{split}
\int_{K_0-B_l} |Tb_{lm}| d\lambda &
\leq \Cr{constant:-2nd-Derivative-Bound-Newtonian-Potential}\rho\int_{K_0 -B_l} \frac{1}{r^{n+1}}\int_{\Upsilon_l} |b_{lm}| d\lambda \\&
= \Cr{constant:-2nd-Derivative-Bound-Newtonian-Potential} \rho \int_{\frac{\delta}{2}}^{\infty}  \frac{1}{r^2} dr\int_{|\omega|=1} d\omega \int_{\Upsilon_l} |b_{lm}|\\&
= \Cr{constant:-2nd-Derivative-Bound-Newtonian-Potential} \rho n \omega_n \frac{2}{\rho} \int_{\Upsilon_l} |b_{lm}| \\&
= 2\Cr{constant:-2nd-Derivative-Bound-Newtonian-Potential} n \omega_n \int_{\Upsilon_l} |b_{lm}|,
\end{split}
\end{align}
where in the third step we used that the volume of the unit sphere in $\R^n$ is $n \omega_n$.
Define $\Cl{constant:Integral-Bound-Calderon-Zygmund} = 2 \Cr{constant:-2nd-Derivative-Bound-Newtonian-Potential}n \omega_n = 2 n^2(n+5)$.
Next write $F^*= \cup B_l$, $G^* = K_0 - F^*$.
For each $l,m \in \mathbb{N}$ we have from \cref{equation:Tblm-estimate-on-K0-minus-Bl}:
\begin{align*}
    \int_{G^*} |Tb_{lm}| \leq \Cr{constant:Integral-Bound-Calderon-Zygmund} \int_{\Upsilon_l} |b_{lm}|.
\end{align*} 
Taking the limit as $m \rightarrow \infty$ we have the following:
\begin{align}
    \label{equation:Tbl-on-Gstar-bound}
    \int_{G^*} |Tb_l| \leq \Cr{constant:Integral-Bound-Calderon-Zygmund} \int_{\Upsilon_l} |b_{l}|. 
\end{align}
We now sum over $l$ and obtain:
\begin{align*}
     \int_{G^*} |Tb| \leq \Cr{constant:Integral-Bound-Calderon-Zygmund} \sum_{l=1}^{\infty} \int_{\Upsilon_l} |b_l|
     =
     \Cr{constant:Integral-Bound-Calderon-Zygmund}\int_{\cup \Upsilon_l} |b| \leq \Cr{constant:Integral-Bound-Calderon-Zygmund} \int_{K_0} |b| \leq \Cr{constant:Integral-Bound-Calderon-Zygmund} \int_{K_0} |f| + |g| \leq \int_{K_0} 2 \Cr{constant:Integral-Bound-Calderon-Zygmund} |f|,
\end{align*}
where in the first step we used $|Tb|=\sum_{l=1}^\infty |Tb_l|$ from \cref{equation:bl-definition} together with \cref{equation:Tbl-on-Gstar-bound}

Applying \cite[Lemma 9.7]{Gilbarg1998} we obtain:
\begin{align}
\label{dist G* bound}
\left| \left\{
x \in G^* \text{ satisfying } |Tb(x)| > \frac{t}{2}
\right\} \right| 
\leq \frac{2}{t} \int_{G^*} |Tb| \leq \frac{4 \Cr{constant:Integral-Bound-Calderon-Zygmund} \|{f}_{L^1}}{t}.
\end{align}
We bound the distribution function of $Tb$ on $F^*$ by simply bounding the measure of $F^*$. 
By \cite[p.234]{Gilbarg1998}:
\begin{align}
\label{balls and cubes}
|F^*| \leq \Cl{constant:-F^*-measure-bound} |\Upsilon|,
\end{align}
where $\Cr{constant:-F^*-measure-bound} = \omega_n n^{n/2}$.
\\
We use \cref{nice property of bad functions} to deduce that for each subcube $\Upsilon_l \in \Upsilon$:
\begin{align*}
    |\Upsilon_l| \leq \frac{\int_{\Upsilon_l} f}{t}.
\end{align*}
Summing over all subcubes in $\Upsilon$ and applying \cref{balls and cubes} gives us:
\begin{align}
\label{dist F* bound}
   |F^*| \leq \Cr{constant:-F^*-measure-bound} \frac{\|{f}_{L^1}}{t}.
\end{align}

The prior bounds \cref{dist F* bound} and \cref{dist G* bound} imply that:
\begin{align*}
    \mu_{Tb}(\frac{t}{2})\leq \frac{(4 \Cr{constant:Integral-Bound-Calderon-Zygmund}+ \Cr{constant:-F^*-measure-bound} )\|{f}_{L^1}}{t}.
\end{align*}
We substitute this into \cref{good+bad} and together with \cref{equation:mu-Tg-estimate} obtain that:
\begin{align*}
   \mu_{Tf}(t) &\leq \mu_{Tg} \left(\frac{t}{2} \right) + \mu_{Tb}\left( \frac{t}{2} \right)
             \leq  (2^{n+2} + 4 \Cr{constant:Integral-Bound-Calderon-Zygmund}+ \Cr{constant:-F^*-measure-bound}) \frac{\|{f}_{L^1}}{t}.
\end{align*}
Thus, we have proven \cref{Distribution Function Bounds} to hold for the following values of $\Cr{constant:mu-t-L2-norm}$ and $\Cr{constant:mu-t-L1-norm}$:
\begin{align}
\label{distribution-function-bounds}
     \Cr{constant:mu-t-L1-norm} = 2^{n+2} + 4 \Cr{constant:Integral-Bound-Calderon-Zygmund}+ \Cr{constant:-F^*-measure-bound},  \Cr{constant:mu-t-L2-norm} =1.
\end{align}
Applying the Marcinkiewicz interpolation theorem finishes the proof for $1 \leq p \leq 2$.
We obtain that:
\begin{align*}
 \|{D_{ij} w} \leq C_{Marcinkiewicz}(1,2,p) T_1^\alpha T_2^{1-\alpha} \|{f}_{L^p}.
\end{align*}
As in \cite[Theorem 9.9]{Gilbarg1998}, this is extended to $p \geq 2$ by duality.
Namely, the Calderon-Zygmund Constant for $p$ is the same as that for the conjugate of $p$, which completes the proof.
\end{proof}
We state the estimate relevant to the rest of this section:

\begin{remark}
\label{corollary:Calderon-Zygmund-3.1}
   We are interested in the constant in the case $p=4$ and $n=3$, we have that:
   \begin{align*}
    \Cl{constant:Calderon-Zygmund-3.1}:= C_{Calderon-Zygmund}(3,4) = 193
   \end{align*}
\end{remark}

We are specifically interested in this result in the case of Compactly Supported Functions.
\begin{corollary}[Corollary 9.10 in \cite{Gilbarg1998}]
\label{theorem:calderon-zygmund}
    Define $\|{D^2 u}_{L^p}:= \|{(|D^2u|_{Frob}) }_{L^p}$
    Let $u \in W^{2,4}_{0}(\Omega)$, $1 < p < \infty$, where $W^{2,p}_0(\Omega)$ denotes the $L^4$-Sobolev space with two weak derivatives and compact support on $\Omega \subset \R^n$.
    Then
    \begin{align*}
        \|{D^2 u}_{L^p}
        \leq
       n^2 C_{\text{Calderon-Zygmund}}(n,p)
        \|{\Delta u}_{L^p},
    \end{align*}
    where $C_{\text{Calderon-Zygmund}}(n,p)$ was defined in \cref{theorem:calderon-zygmund-potential-formulation}.
\end{corollary}
\begin{proof}
We have that $u$ is the Newtonian potential of $\Delta u$, because $u$ has compact support (by \cite[Equation 2.17]{Gilbarg1998}), so we can apply \cref{theorem:calderon-zygmund-potential-formulation}.
\begin{align}
\begin{split}
\label{Frobenius Calderon-Zygmund}
  \|{D^2u}_{L^p}:= \|{\sqrt{\sum_{i=1, j=1}^n(D_{ij}u)^2}}_{L^p} &\leq \sqrt{\|{\sum_{i=1. j=1}^{n} |D_{ij}(u)|^2}}_{L^p}
  \\&\leq \sqrt{\left(\sum_{i=1, j=1}^n \|{D_{ij}(u)}_{L^p}\right)^2}
  \\&\leq n^2 C_{\text{Calderon-Zygmund}} (n,p) \|{\Delta u}_{L^p}
  \end{split}
\end{align}
In the second inequality we apply \cite[Theorem 202]{HardyLittlewoodPolya1934}
\end{proof}

\begin{corollary}
    \label{corollary:D+D^2-regularity-estimate}
    Let $u \in W^{2,p}_0(\Omega)$, $1 < p < \infty$.
    Then
    \begin{align*}
        \|{D^2 u}_{L^p}
        +
        \|{D u}_{L^p}
        \leq
        \Cl{D+D^2-regularity-estimate}
        \|{\Delta u}_{L^p},
    \end{align*}    
    where $\Cr{D+D^2-regularity-estimate}= n^2C_{Calderon-Zygmund(n,p)}(n\Cr{Poincaré Inequality}+1)$.
    
\end{corollary}

\begin{proof}
    The Calderon-Zygmund estimate, \cref{theorem:calderon-zygmund}, bounds $\|{D^2u}_{L^p}$, it remains to bound $\|{Du}_{L^p}$.
    We have that $\frac{\partial}{\partial x_i} u$ is compactly supported for $i \in \{1,\dots,n\}$.
    Thus, by the Poincaré inequality \cref{Poincaré Inequality}:
    \begin{align}
    \label{Partial Derivative Poincaré Estimate}
        \|{D_i u}_{L^p}
        \leq
        \Cl{Poincaré Inequality}
        \|{D D_i u}_{L^p},
    \end{align}

Where $\Cr{Poincaré Inequality} = (\frac{1}{\omega_n} |\Omega|)^\frac{1}{n}$.
Hence
\begin{align*}
\|{Du}_{L^p}
= 
\|{ \left( \sum_{i=1}^n |D_i(u)|^2 \right)^\frac{1}{2}}_{L^p}
\leq
\sqrt{\|{\sum_{i=1}^n |D_i(u)|^2}_{L^p}} 
\leq 
\sqrt{\left(\sum_{i=1}^n \|{D_i(u)}_{L^p}\right)^2} 
\leq n \max_{1\leq i \leq n} \left(\|{D_i(u)}_{L^p} \right) .
\end{align*}
In the second inequality we apply \cite[Theorem 202]{HardyLittlewoodPolya1934} inside the square root. 
We combine the above two estimates to conclude:
\begin{align*}
    \|{D^2u}_{L^p} + \|{Du}_{Lp} 
    &\leq (n \Cr{Poincaré Inequality} +1) \|{D^2u}_{L^p}
     \leq n^2 C_{Calderon-Zygmund}(n,p) (n\Cr{Poincaré Inequality} +1) \|{\Delta u}_{L^p}.
    \qedhere
\end{align*}
\end{proof}
The constant that we have derived in \cref{corollary:D+D^2-regularity-estimate} involves a constant term  $\Cr{Poincaré Inequality}$ which is dependent on the volume of $\Omega$, and in our case we will be applying \cref{corollary:D+D^2-regularity-estimate} to $\Omega = \tilde{Q} = [-1,2]^3$, this will be motivated shortly.
\\
Denote $\Cl{Poincaré-Inequality-Q-tilde} = (\frac{81}{4\pi})^{1/3}$, as the value of $\Cr{Poincaré Inequality}$ with $\Omega = \tilde{Q}$

\subsection{Estimates on 3-Torus with the flat metric}
\label{subsection:estimates-on-t3-with-flat-metric}
\begin{theorem}
    \label{theorem:t3-schauder}
    From now we specialise to the case $p=4$
    Let $u \in W^{2,4}(T^3)$.
    Then
    \begin{align}
        \label{equation:schauder-estimate-on-T3}
        \|{u}_{W^{2,4}(T^3)}
        \leq
        \Cl{constant:Schauder-Estimate}
        \left(
        \|{\Delta u}_{L^4(T^3)}
        +
        \|{u}_{L^4(T^3)}
        \right),
    \end{align}.
\end{theorem}
Where $\Cr{constant:Schauder-Estimate}$ is defined in \cref{Explicit Constant for Schauder Estimate} 
\begin{proof}
    As usual we view $u \in W^{2,4}(T^3)$ as a function in $W^{2,4}(\R^3)$ that is periodic in each of the coordinate directions. 
    Using a cutoff function we can apply the results from the previous section which requires "zero boundary".  
    Let $Q:=[0,1]^3 \subset \R^3$ and $\tilde{Q}:=[-1,2]^3 \subset \R^3$.
    Let $\chi: \tilde{Q} \rightarrow [0,1]$ be twice differentiable such that $\chi(x)=1$ for all $x \in Q$ and $\chi(y)=0$ for all $y \in \partial \tilde{Q}$. We define such a $\chi$ explicitly in \cref{proposition:cut-off-estimates}.
    Then
    \begin{align}
        \label{eqn:elliptic-regularity-torus-proof-step}
        \begin{split}
        \|{u}_{W^{2,4}(T^3)}
        &=
        \|{\chi u}_{W^{2,4}(Q)},
        \\
        &\leq
        \|{\chi u}_{W^{2,4}(\tilde{Q})},
        \\
        &=
        \|{\chi u}_{L^4(\tilde{Q})}
        +
        \|{D(\chi u)}_{L^4(\tilde{Q})}
        +
        \|{D^2(\chi u)}_{L^4(\tilde{Q})},
        \\
        &\leq
        \|{\chi u}_{L^4(\tilde{Q})}
        +
        \Cr{D+D^2-regularity-estimate}\|{\Delta (\chi u)}_{L^4(\tilde{Q})}.
        \end{split}
    \end{align}
    where we used \cref{corollary:D+D^2-regularity-estimate} in the last step.
    Here we find for the last term:
    \begin{align}
        \label{eqn:delta(chi.u)-estimate}
        \|{\Delta (\chi u)}_{L^4(\tilde{Q})}
        &\leq
        \underbrace{
        \|{\Delta \chi}_{L^\infty(\tilde{Q})} 
        }_{\leq \Cl{delta-chi-bound}}
        \|{u}_{L^4(\tilde{Q})}
        +
        2\|{D\chi \cdot Du}_{L^4(\tilde{Q})}
        +
        \underbrace{
        \|{\chi}_{L^\infty(\tilde{Q})}
        }_{\leq 1}
        \|{\Delta u}_{L^4(\tilde{Q})},
    \end{align}
    where $\chi$ is the cut-off function and constants $\Cr{delta-chi-bound}$, $\Cl{D-chi-bound}$ and $\Cl{D^2-chi-bound}$ are defined in \cref{proposition:cut-off-estimates}. We now bound the middle term:
    \begin{align}
    \label{Middle Term Estimate}
    \begin{split}
        \|{(D\chi) \cdot (D u)}_{L^4(\tilde{Q})}
        &\leq
        27^\frac{1}{12}
        \|{(D\chi) \cdot (Du) }_{L^6(\tilde{Q})}\\
        &\leq 
        27^\frac{1}{12}\Cr{sobolev-embedding}\|{D((D\chi) \cdot (D u))}_{L^2(\tilde{Q})}
        \\
        &\leq
        27^\frac{1}{12}
        \Cr{sobolev-embedding}
        \left(
        \|{D^2 \chi}_{C^0(\tilde{Q})} \cdot \|{D u}_{L^2(\tilde{Q})}
        +
        \|{D \chi}_{C^0(\tilde{Q})} \cdot \|{D^2 u}_{L^2(\tilde{Q})}
        \right)
        \\
        &\leq
        27^\frac{1}{12}
        \max(\|{D^2\chi}_{L^\infty(\tilde{Q})}, \|{D\chi}_{L^\infty(\tilde{Q})})\Cr{sobolev-embedding} \|{u}_{W^{2,2}(\tilde{Q})}
        \\
        &\leq
        27^\frac{1}{12}
        \Cr{D^2-chi-bound} \Cr{sobolev-embedding} \|{u}_{W^{2,2}(\tilde{Q})}.
    \end{split}
    \end{align}
    In the first step we use Hölder's inequality;
    in the second step we use the Sobolev Embedding Theorem, \cref{Sobolev Embedding};
    third we use the product rule and last we use \cref{proposition:cut-off-estimates}.
    Since it arises fairly frequently denote 
    \begin{align*}
    \Cl{Holder Q-Qtilde}= 27^\frac{1}{12}.
    \end{align*}
In order to prove \cref{theorem:t3-schauder}, we require the following lemma:
    \begin{lemma}[Lemma 8.2.3 in \cite{Jost2014}]
    \label{Du estimate on L^2}
    Let $u$ be a weak solution of $\Delta u=f$ with $f \in L^2(\Omega)$.
    We then have for any $\Omega' \subset \Omega$ whose closure is contained in $\Omega$
    \begin{align*}
        \|{Du}_{L^2(\Omega')} \leq \frac{\sqrt{17}}{\delta^2}\|{u}_{L^2(\Omega)} + \delta^2\|{\Delta u}_{L^2(\Omega)},
    \end{align*} 
    where $\delta = d(\Omega', \partial\Omega)$
    \end{lemma}
    
    \begin{remark}
    \label{remarkable remark}
    In the following we use \cref{Du estimate on L^2} with $\Omega' = Q, \Omega = \tilde{Q}, \delta= 1$.
    \end{remark}

We return to the proof of \cref{equation:schauder-estimate-on-T3}.
We perform the tedious computation of combining the above lemma with \cref{eqn:elliptic-regularity-torus-proof-step} and \cref{eqn:delta(chi.u)-estimate}
\begin{align*}
\|{u}_{L^4(T^3)} &\leq \|{\chi u}_{L^4(\tilde{Q})} + \Cr{D+D^2-regularity-estimate} \|{\Delta(\chi u)}_{L^p(\tilde{Q})}
\\&
\leq  27\|{u}_{L^4(Q)} 
+ \Cr{D+D^2-regularity-estimate}(\Cr{delta-chi-bound} \|{u}_{L^4(\tilde{Q})} +  2\Cr{Holder Q-Qtilde}\Cr{D^2-chi-bound} \Cr{sobolev-embedding} \|{u}_{L^2_2(\tilde{Q})} + \|{\Delta u}_{L^4(\tilde{Q})}) \\&
\leq  27(1+ \Cr{D+D^2-regularity-estimate}\Cr{delta-chi-bound}) \|{u}_{L^4(Q)} + 27 \Cr{D+D^2-regularity-estimate} \|{\Delta u}_{L^4(Q)} + 54\Cr{D+D^2-regularity-estimate} \Cr{Holder Q-Qtilde}\Cr{D^2-chi-bound}\Cr{sobolev-embedding} \|{u}_{W^{2,2}(Q)} \\&
\leq  27(1+\Cr{D+D^2-regularity-estimate}\Cr{delta-chi-bound}) \|{u}_{L^p(Q)} + 27\Cr{D+D^2-regularity-estimate}\|{\Delta u}_{L^4(Q)}+ 54\Cr{D+D^2-regularity-estimate} \Cr{Holder Q-Qtilde}\Cr{D^2-chi-bound}\Cr{sobolev-embedding}(\|{u}_{L^2(T^3)} + \|{Du}_{L^2(T^3)}+\|{D^2u}_{L^2(Q)}) \\&
\leq  27(1+\Cr{D+D^2-regularity-estimate}\Cr{delta-chi-bound} +2\Cr{D+D^2-regularity-estimate}\Cr{Holder Q-Qtilde}\Cr{D^2-chi-bound}\Cr{sobolev-embedding}) \|{u}_{L^4(Q)} + 27(\Cr{D+D^2-regularity-estimate} + 2\Cr{D+D^2-regularity-estimate}\Cr{Holder Q-Qtilde}\Cr{D^2-chi-bound}\Cr{sobolev-embedding})\|{\Delta u}_{L^4(Q)} + 54\Cr{D+D^2-regularity-estimate}\Cr{Holder Q-Qtilde}\Cr{D^2-chi-bound}\Cr{sobolev-embedding} \|{Du}_{L^2(Q)}\\&
\leq 27(1+\Cr{D+D^2-regularity-estimate}\Cr{delta-chi-bound} +2\Cr{D+D^2-regularity-estimate}\Cr{Holder Q-Qtilde}\Cr{D^2-chi-bound}\Cr{sobolev-embedding}) \|{u}_{L^4(Q)} + 27(\Cr{D+D^2-regularity-estimate} + 2\Cr{D+D^2-regularity-estimate}\Cr{Holder Q-Qtilde}\Cr{D^2-chi-bound}\Cr{sobolev-embedding})\|{\Delta u}_{L^4(Q)}  \\&\quad+54\Cr{D+D^2-regularity-estimate}\Cr{Holder Q-Qtilde}\Cr{D^2-chi-bound}\Cr{sobolev-embedding}(\sqrt{17}\|{u}_{L^2(\tilde{Q})} + \|{\Delta u }_{L^2(\tilde{Q})})\\&
\leq \{27(1+\Cr{D+D^2-regularity-estimate}\Cr{delta-chi-bound} + (2+54(27)\sqrt{17}) \Cr{D+D^2-regularity-estimate}\Cr{Holder Q-Qtilde} \Cr{D^2-chi-bound} \Cr{sobolev-embedding} \} \|{u}_{L^4(Q)} + \{27(\Cr{D+D^2-regularity-estimate} + 2+ 54(27) \Cr{D+D^2-regularity-estimate} \Cr{Holder Q-Qtilde} \Cr{D^2-chi-bound} \Cr{sobolev-embedding} \} \|{\Delta u}_{L^4(Q)}
\end{align*}
In the first step we recall \cref{eqn:elliptic-regularity-torus-proof-step}, second we use \cref{eqn:delta(chi.u)-estimate} and substitute \cref{Middle Term Estimate}, third we use the periodicity assumption, fifth we use the fact that $C_{Calderon-Zygmund}(n,2) =1$ (when the matrix norm is the Frobenius Norm, see \cref{corollary:Calderon-Zygmund-3.1}) for any choice of n and Holder's inequality, in the sixth step we use \cref{Du estimate on L^2} in accordance with \cref{remarkable remark} and the last step is similar to the third step. This concludes the proof of \cref{equation:schauder-estimate-on-T3}, with
\begin{align}
\label{Explicit Constant for Schauder Estimate}
\Cr{constant:Schauder-Estimate} &=  27(1+\Cr{D+D^2-regularity-estimate}\Cr{delta-chi-bound} + (2+54(27)\sqrt{17}) \Cr{D+D^2-regularity-estimate}\Cr{Holder Q-Qtilde} \Cr{D^2-chi-bound} \Cr{sobolev-embedding} .
\end{align}
\end{proof}

We have shown \cref{theorem:t3-schauder}, our goal is to prove \cref{proposition:injectivity-estimate-flat-metric} and the only difference is the $\|{u}_{L^p(T^3)}$ term on the right hand side.
In the rest of the section we explain how to bound this term by $\Delta u$ as well, this is reliant on the assumption in \cref{proposition:injectivity-estimate-flat-metric} that $f$ is $L^2$-Orthogonal to $\Ker \Delta$\\
The following result is well known:

\begin{lemma}
\label{orthonormal-eigenbasis-Laplace-Hilbert}
Assume that $u \in L^2(T^3)$, with $\int_{T^3} u =0$, then:
\begin{align*}
\|{u}_{L^2(T^3)} \leq \frac{1}{4\pi^2}\|{\Delta u}_{L^2(T^3)}.
\end{align*}
\end{lemma}
\begin{proof}
By \cite[p.30]{Chavel1984}, the smallest eigenvalue of $\Delta$ acting on functions with mean zero is $4 \pi^2$.
Writing $u$ in an orthonormal eigenbasis for $\Delta$ proves the claim.
\end{proof}

In the following we appeal to \cref{other-sobolev-embedding}:

We now have gathered enough tools to prove the main result of the section \cref{proposition:injectivity-estimate-flat-metric}.

\begin{proof}[Proof of \cref{proposition:injectivity-estimate-flat-metric}]
We have
\begin{align*}
\begin{split}
\|{u}_{W^{2,4}(T^3)} &\leq \Cr{constant:Schauder-Estimate}( \|{u}_{L^4(T^3)} + \|{\Delta u}_{L^4(T^3)})\\&
\leq \Cr{constant:Schauder-Estimate} \Cr{Sobolev-Embedding-Cube} \|{u}_{L^2(T^3)} +\Cr{constant:Schauder-Estimate} \Cr{Sobolev-Embedding-Cube} \|{Du}_{L^2(T^3)} + \Cr{constant:Schauder-Estimate} \|{\Delta u}_{L^4(T^3)}\\&
\leq \Cr{constant:Schauder-Estimate} \Cr{Sobolev-Embedding-Cube} \|{u}_{L^2(T^3)} + \Cr{constant:Schauder-Estimate} \|{\Delta u}_{L^4(T^3)} +\Cr{constant:Schauder-Estimate} \Cr{Sobolev-Embedding-Cube}(\sqrt{17} \|{u}_{L^2(\tilde{Q})} + \|{\Delta u}_{L^2(\tilde{Q}})\\&
\leq  \Cr{constant:Schauder-Estimate} \Cr{Sobolev-Embedding-Cube}(1 +27\sqrt{17} ) \|{u}_{L^2(T^3)} + \Cr{constant:Schauder-Estimate}(1+27\Cr{Sobolev-Embedding-Cube}) \|{\Delta u}_{L^4(T^3)} \\&
\leq \left[
\frac{1}{4\pi^2}\{\Cr{constant:Schauder-Estimate} \Cr{Sobolev-Embedding-Cube}(1 + 27\sqrt{17} )\} + \Cr{constant:Schauder-Estimate}(1+27\Cr{Sobolev-Embedding-Cube})
\right] \|{\Delta u}_{L^4(T^3)}.
\end{split}
\end{align*}
In the first step we appeal to \cref{equation:schauder-estimate-on-T3}, we next use \cref{other-sobolev-embedding}, in the third step we use \cref{Du estimate on L^2} with \cref{remarkable remark}, in the last step we use \cref{orthonormal-eigenbasis-Laplace-Hilbert} together with $\|{\Delta u}_{L^2} \leq \vol(T^3)^{\frac{p-2}{2p}}\|{\Delta u}_{L^p}$ and $\vol(T^3)=1$.
This proves the proposition with 
\begin{align}
    \label{equation:flat-injectivity-constant}
    \Cr{const:injectivity-main-flat} = \frac{1}{4\pi^2}\{\Cr{constant:Schauder-Estimate} \Cr{Sobolev-Embedding-Cube}(1 + 27\sqrt{17} )\} + \Cr{constant:Schauder-Estimate}(1+27\Cr{Sobolev-Embedding-Cube}\}).
\end{align}
\end{proof}

\subsection{Estimates on 3-Torus with respect to a non-flat metric}
\label{subsection:estimates-on-t3-with-non-flat-metric}

We now turn to proving the injectivity estimate on $T^3$ with respect to non-flat metrics, namely \cref{theorem:T3-injectivity-estimate}.\\
Throughout the section, let $\gflat$ denote the flat metric on $T^3$ and let $g$ denote another metric.\\
In this section, we work exclusively with the manifold $T^3$. Thus we omit it from the notation introduced to describe norms on Sobolev and $L^p$ spaces,
For example:
\begin{align}
\text{We write }  g, W^{2,p} \text{ instead of }  W^{2,p}(T^3.g)
\end{align}
If the metric is not specified in a calculation assume that it is $\gflat$. 
In calculations we use the coordinates $x_1,x_2,x_3$ which have the property that $e_i:=\frac{\partial}{\partial x_i}$ form an orthonormal basis with respect to $\gflat$, but not necessarily with respect to $g$.
We write $g_{ij}=g \left( \frac{\partial}{\partial x_i}, \frac{\partial}{\partial x_j} \right)$ for the matrix representation of $g$ in this coordinate basis.

The claim \cref{theorem:T3-injectivity-estimate} will be established by showing the following intermediate inequalities, which illustrates the rough strategy of the proof.
\begin{align}
\label{overall strategy for injective estimate with arbitrary-ish riemannian metric}
\|{u}_{g, W^{2,4}}\leq \Cl{Step 1 3.3} \|{u}_{\gflat, W^{2,4}}\leq \Cr{const:injectivity-main-flat}\Cr{Step 1 3.3} \|{\Delta u}_{\gflat, L^4}\leq \Cl{Step 3 3.3} \|{\Delta u}_{g, L^4} + C(\delta) \|{u}_{\gflat, W^{2,4}}
\end{align}
In the last inequality if we take $\delta>0$ to be sufficiently small can absorb the remainder term to obtain our desired inequality.
To show our application in \cref{section:application-nowhere-vanishing-1-forms}, it also turns out to be necessary to find an explicit constant for the following inequality:
\begin{align}
\|{u}_{g, W^{2,4}} \geq \Cl{Step 1 3.3 Lower Bound} \|{u}_{\gflat, W^{2,4}}
\end{align}
Before we proceed we provide some linear algebra estimates that will be useful later.
\begin{definition}
Define the following norm on matrx valued functions on $T^3$
\begin{align*}
    \|{A}_{C^0} &= \max_{1\leq i,j \leq 3, x \in T^3} |A_{ij}(x)|\\
    \|{A}_{C^1} &= \|{A}_{C_0} + \|{dA}_{C^0}
\end{align*}
Note this norm is not submultiplicative.
\end{definition}
\begin{lemma}
\label{Determinant Estimate}
If $\|{g-\gflat}_{C^0} \leq \delta$, then we have the following bound for the determinant:
\begin{align*}
\Cl{Lower Determinant Estimate} \leq  \det(g) \leq \Cl{Upper Determinant Estimate},
\end{align*}
where $\Cr{Lower Determinant Estimate}$ and $ \Cr{Upper Determinant Estimate}$ are defined in \cref{equation:det-constants}.
\end{lemma}
\begin{proof}
The assumption implies
$1- \delta \leq g_{ii} \leq 1+\delta,
|g_{ij}| \leq \delta$ for $i \neq j$.
The claim then easily follows from the explicit formula for the determinant of a $3 \times 3$ matrix with 
\begin{align}
\label{equation:det-constants}
\Cr{Upper Determinant Estimate} = (1+\delta)^3 + 2\delta^3 + 3(1+\delta) \delta^2, \Cr{Lower Determinant Estimate} = (1-\delta)^3 - 2\delta^3 - 3(1+\delta)\delta^2.
\end{align}
\end{proof}

The following two lemmas involving the inverse of the Riemmanian metric also play a role:

\begin{lemma}
Suppose $\|{g-\gflat}_{C^1} \leq \delta$ and that $\delta < \frac{1}{6}$, then we have the following:
\begin{align}
\label{inverse metric bound}
\|{g^{-1}-\gflat}_{C^1} \leq 2 \delta.
\end{align}
\begin{proof}
We prove the $C^0$ bound and $C^1$ bound seperately. 
Denote $g= I - A$, we intend to use the series expansion of $g^{-1} = (I - A)^{-1}$ and then use submultiplicativity, the $C^0$ and $C^1$ norms on matrices that we defined are not submultiplicative so we compare them to different matrix norms. 
\begin{align*}
    \|{(I-A)^{-1} -I}_{C^0} \leq \|{(I-A)^{-1} -I }_{1}:= \sup_{j \in \{1,2,3\}} \sum_{i=1}^3 |A_{ij}| 
\end{align*}
\begin{align*}
\|{(I - A)^{-1} - I}_{C^0} \leq \|{\sum_{i=0}^\infty A^i-I}_{1} \leq \sum_{i=1}^\infty (3\delta) ^ i \leq 6 \delta,
\end{align*}
We now bound the derivative:
We first set $B:= I-A$ and we differentiate the identity $BB^{-1} =I$ to obtain:
\begin{align}
B'(B^{-1}) + B(B^{-1})' = 0 .
\end{align}
Rearranging this identity and multiplying both sides by $B^-1$ to the left gives the following:
\begin{align}
(B^{-1})' = -B^{-1}B'B^{-1} .
\end{align}
Using the submultiplicative property of the 1-norm:
\begin{align}
\|{(B^{-1})'}_{1} \leq (\|{B^{-1}}_{1})^2 \|{B'}_{1}.
\end{align}
Using our prior work we have that:
\begin{align}
\|{(B^{-1})'}_{C^0} \leq \|{(B^{-1})'}_{1} \leq 36 \delta^3. 
\end{align}
We obtain that:
\begin{align}
\|{I-A}_{C^1} \leq 6 \delta + 36 \delta^3,
\end{align}
which gives \cref{inverse metric bound}.
This concludes the proof.
\end{proof}
\end{lemma}
The next lemma compares the pointwise norm of covectors with respect to the flat metric and the perturbed metric. 
\begin{lemma}
\label{comparison lemma vectors}
Assume that V is a vector space equipped with two inner products $g,h$. 
For an \textit{h}-orthonormal basis $e_1, e_2, ..., e_n$ write $g_{ij}= g(e_i,e_j)$
If $|g_{ij} - \delta_{ij}| \leq \delta$, for all $i,j \in \{1,2,3,..., n\}$ then for each $v \in V$ we have that:
\begin{align}
(1-n\delta)^{\frac{1}{2}} |v|_{h}  \leq |v|_g \leq (1+n\delta)^{\frac{1}{2}} |v|_{h}.
\end{align}
\end{lemma}

\begin{proof}
Fix $v \in V$ it follows that $v = \sum_{i=1}^n v_i e_i$
\begin{align}
  \begin{split}
  |g(v,v)-h(v,v)| &= 
  \left|
  \sum_{i=1, j=1}^n v_i v_j [g(e_i, e_j) - h(e_i, e_j)]
  \right|
            \\&\leq \delta  \sum_{i=1, j=1}^n |v_i v_j| 
            \\& \leq \delta (\sum_{i=1}^n |v_i|)^2 
            \\& \leq \delta n h(v,v) 
  \end{split}
\end{align}
In the third step we observe that $|v_i v_j| = |v_i| |v_j|$ and manipulate the sum accordingly and in the final step we apply norm equivalence between the 1-norm and the 2-norm.
By using the triangle inequality we obtain that:
\begin{align}
(1-n\delta) h(v,v) \leq g(v,v) \leq (1+n\delta) h(v,v).
\end{align}
Taking the square root in each term of the above gives the required result. 
\end{proof}
Slight adjustments to the above lemma gives us the analogous claim for tensors.
\begin{lemma}
\label{comparison lemma covectors}
Suppose $w \in T_x^*(T^3)$, where $T_x^*(T^3)$ is the cotangent space of $T^3$ at $x$ and $|g^{ij}(x)-\delta_{ij}| < \delta$, we obtain the following bound on pointwise norms.
\begin{align*}
\Cl{Covector Lower Estimate} |w|_{\gflat} \leq|w|_g \leq \Cl{Covector Estimate} |w|_{\gflat},
\end{align*}
with $\Cr{Covector Estimate} = (1+3\delta)^{\frac{1}{2}}$ and $\Cr{Covector Lower Estimate} = (1-3\delta)^{\frac{1}{2}}$.
\end{lemma}

\begin{proof}
$T_x^*(T^3)$ is a vector space with orthonormal basis $dx_i$ for $i \in \{1,2,3\}$.
Applying \cref{comparison lemma vectors} gives the claim.  
\end{proof}
We now prove a similar claim to \cref{comparison lemma covectors} but for covariant 2 tensors.


\begin{lemma}
\label{lemma:2-tensor-metric-comparison}
Suppose $ T = \sum_{i=1, j=1}^3 T_{ij} dx_i \otimes dx_j \in (T_x^*T^3)^{\otimes 2} $ and $|g^{ij}(x)-\delta_{ij}| < \delta$, then we have the following bound on pointwise norms.
\begin{align*}
\Cl{Co2Tensor Lower Estimate}|T|_{\gflat}\leq |T|_{g} \leq \Cl{Co2Tensor Estimate}|T|_{\gflat}.
\end{align*}
\end{lemma}
\begin{proof}
We would like to use \cref{comparison lemma vectors}.
Note that $dx_i \otimes dx_j$ with $i,j \in \{1,2,3\}$ forms an orthonormal basis with respect to $\gflat$.
It suffices to establish a bound for $|g^{ij,kl}- \delta_{ij,kl}|$,
where we used the notation
$g^{ij,kl}=\langle dx_i \otimes dx_j, d_k \otimes dx_l \rangle = g^{ik} g^{jl}$.
\begin{align}
\begin{split}
|g^{ij}g^{kl} - \delta_{ij}\delta_{kl}| &\leq |g^{ij}g^{kl} - g^{ij}\delta_{kl}+ g^{ij} \delta_{kl} - \delta_{ij}\delta_{kl}| \\&\leq |g^{ij}g^{kl} - g^{ij}\delta_{kl}| +  |g^{ij}\delta_{kl} - \delta_{ij} \delta_{kl}| \\&\leq |g^{ij}||g^{kl}- \delta_{kl}| + |\delta_{kl}||g^{ij} - \delta_{ij}| \\&\leq (2+2\delta) 2\delta .
\end{split}
\end{align}
Now applying \cref{comparison lemma vectors} yields the claim with:
\begin{align}
\Cr{Co2Tensor Lower Estimate} = \sqrt{1- 9(2+2\delta)2\delta},\\
\Cr{Co2Tensor Estimate} = \sqrt{1+ 9(2+2\delta)2\delta}.
\end{align}
\end{proof}

We also require a bound for the Christoffel Symbols:
\begin{lemma}
\label{lemma:christoffel-bound}
Suppose $\|{g - \gflat}_{C^1} \leq \delta$. Then for $i,j,m \in \{1,2,...,n\}$ we have the following:
\begin{align}
\label{Christoffel Symbol Estimate}
|\Gamma^m_{ij}| \leq \Cl{Christoffel Bound} = 3^2 \delta^2.
\end{align}
\end{lemma}
\begin{proof}
We use the usual formula for the Christoffel Symbols of the Levi-Civita Connection:
\begin{align}
\label{Christoffel Symbols}
\Gamma^k_{ij} = \sum_{n=1}^3 \frac{1}{2} g^{kn} \left(
\frac{\partial g_{ni}}{\partial x_j} + \frac{\partial g_{nj}}{\partial x_i} - \frac{\partial g_{ij}}{\partial x_n}
\right).
\end{align}
By the Triangle Inequality we arrive at the following:
\begin{align}
\begin{split}
|\Gamma^k_{ij}| &\leq \sum_{n=1}^3 \frac{1}{2} |g^{kn}| \left|
\frac{\partial g_{ni}}{\partial x_j} + \frac{\partial g_{nj}}{\partial x_i} - \frac{\partial g_{ij}}{\partial x_n}\right|
\leq
\sum_{n=1}^3 \frac{1}{2} (2\delta)(3\delta)
=
3^2 \delta^2
\end{split}
\end{align}
We use \cref{inverse metric bound} and our assumption in the second step. 
\end{proof}

\begin{theorem}
\label{theorem:L^p_2-norm-comparison}
With the assumptions in \cref{theorem:T3-injectivity-estimate}, we have the following:
\begin{align}
 \label{Step 1: T3 Injectivity Estimate}
    \Cr{Step 1 3.3 Lower Bound} \|{u}_{\gflat, W^{2,p}}\leq \|{u}_{g,W^{2,p}} \leq \Cr{Step 1 3.3} \|{u}_{\gflat, W^{2,p}}.
\end{align}
\end{theorem}

Where:
\begin{align}
\Cr{Step 1 3.3} = \left( \Cr{S1P1 T3 Injectivity Estimate} (1+3\sqrt{3}\Cr{Christoffel Bound}) \Cr{Co2Tensor Estimate} \Cr{Upper Determinant Estimate}^{\frac{1}{2p}} \right),
\Cr{Step 1 3.3 Lower Bound}= \Cr{S1P1 T3 Injectivity Estimate Lower Bound} + (1-3\sqrt{3}\Cr{Christoffel Bound}) \Cr{Lower Determinant Estimate}^{\frac{1}{2p}} \Cr{Co2Tensor Lower Estimate}.
\end{align}
\begin{proof}
\textbf{Step 1:}
We begin by proving
\begin{align}
\label{Step 1 Part 1: T3 Injectivity Estimate}
 \Cl{S1P1 T3 Injectivity Estimate Lower Bound} \|{u}_{\gflat, W^{1,p}}\leq \|{u}_{g, W^{1,p}}\leq \Cl{S1P1 T3 Injectivity Estimate}\|{u}_{\gflat, W^{1,p}},
\end{align}
where $\Cr{S1P1 T3 Injectivity Estimate Lower Bound},\Cr{S1P1 T3 Injectivity Estimate}$ are defined in \cref{equation:L^p_1-comparison-constants}.

We first exhibit the upper bound:
\begin{align*}
\begin{split}
\|{u}_{g, W^{1,p}} 
&= 
\left(
\int_{T^3} |u|^p \sqrt{det(g)}
\right)^\frac{1}{p} 
+ 
\left(
\int_{T^3} |du|^p_g \sqrt{det(g)}
\right)^\frac{1}{p} 
\\
&\leq \Cr{Upper Determinant Estimate}^\frac{1}{2p}
\left[
\left(
\int_{T^3} |u|^p
\right)^\frac{1}{p}+
\left(\int_{T^3} |du|^p_g
\right)^\frac{1}{p}
\right]
\\
&\leq \Cr{Covector Estimate} \Cr{Upper Determinant Estimate}^\frac{1}{2p} \left[ \left(
\int_{T^3}|u|^p
\right)^\frac{1}{p}+
\left(\int_{T^3} |du|^p_{\gflat} \right)^\frac{1}{p}
\right],
\end{split}
\end{align*}
where we used \cref{Determinant Estimate} in the second step, and we used \cref{comparison lemma covectors} in the third step. 

The lower bound follows similarly:
\begin{align*}
\begin{split}
\|{u}_{g, W^{1,p}} &= \left(\int_{T^3} |u|^p \sqrt{det(g)}\right)^{\frac{1}{p}} + \left(\int_{T^3} |du|^p_{g} \sqrt{det(g)} \right)^{\frac{1}{p}}
\\&\geq
\Cr{Lower Determinant Estimate}^{\frac{1}{2p}} \left[\left(\int_{T^3} |u|^p\right)^{\frac{1}{p}}+ \left(\int_{T^3} |du|^p_g\right)^{\frac{1}{p}} \right]
\\&\geq
\Cr{Covector Lower Estimate}\Cr{Lower Determinant Estimate}^{\frac{1}{2p}} \left[\left(\int_{T^3} |u|^p\right)^{\frac{1}{p}}+ \left(\int_{T^3} |du|^p_{\gflat}\right)^{\frac{1}{p}} \right],
\end{split}
\end{align*}
where 
\begin{align}
\label{equation:L^p_1-comparison-constants}
\Cr{S1P1 T3 Injectivity Estimate}= \Cr{Covector Estimate} \Cr{Upper Determinant Estimate}^\frac{1}{2p} 
\text{ and }
\Cr{S1P1 T3 Injectivity Estimate Lower Bound} = \Cr{Covector Lower Estimate}\Cr{Lower Determinant Estimate}^{\frac{1}{2p}}.
\end{align}

\textbf{Step 2:}
Pointwise estimate for $\nabla^g \omega$, where $\omega \in \Omega^1(T^3)$.

In what follows, we write $\nabla$ for the Levi-Civita connection of $\gflat$ and $\nabla^g$ for the Levi-Civita connection of $g$.
We compute at a fixed point that we omit from the notation.
We have
\begin{align}
\begin{split}
\label{nabla estimate}
|\nabla^g \omega|_{g} &\leq |\nabla^g \omega - \nabla \omega|_{g} + |\nabla \omega|_{g}
\leq
\Cr{Co2Tensor Estimate}(|\nabla^g\omega - \nabla \omega|_{\gflat} + |\nabla \omega|_{\gflat}).
\end{split}
\end{align}
where we used the triangle inequality and \cref{lemma:2-tensor-metric-comparison}.



By definition
\begin{align}
\label{nabla pointwise estimate}
|\nabla^g \omega - \nabla \omega|_{\gflat}^2 = \sum_{i=1, j=1}^3 |\nabla^g \omega (e_i, e_j) - \nabla \omega (e_i,e_j) |^2,
\end{align}
where we wrote $e_i = \frac{\partial}{\partial x_i}$, which is an orthonormal basis with respect to $\gflat$.
 
Notice we have that:
\begin{align}
\label{simplify nabla estimate}
\begin{split}
\nabla^g \omega (e_i, e_j) - \nabla \omega (e_i, e_j) 
& =
\omega(\nabla_{e_i} e_j) - \omega (\nabla^g_{e_i} e_j) 
=
- \omega (\nabla^g_{e_i} e_j),
\end{split}
\end{align}
where we used that $\nabla_{e_i} e_j=0$ for $i,j \in \{1,2,3\}$.

We will first compute the required estimate for $dx_m$, the standard orthonormal basis for 1-forms on $T^3$, which are dual to the $e_m$ from above.
The same result for $\omega$ will follow immediately. Substituting \cref{Christoffel Symbols} and \cref{simplify nabla estimate} into \cref{nabla pointwise estimate} gives us:
\begin{align}
\label{nabla difference bound}
\begin{split}
|\nabla^g dx_m - \nabla dx_m|^2_{\gflat} &= 
\sum_{i=1, j=1}^3 | dx_m( \nabla^g_{e_i} e_j)|^2 \\& 
=
\sum_{i=1, j=1}^3 
\left| dx_m
\left( 
\sum_{k=1}^3 
\left(
\sum_{n=1}^3  \frac{1}{2} g^{kn} 
\left(
\frac{\partial g_{ni}}{\partial x_j} + \frac{\partial g_{nj}}{\partial x_i} - \frac{\partial g_{ij}}{\partial x_n}
\right)
\right)
e_k
\right)
\right|^2 
\\& 
=
\sum_{i=1,j=1}^3 
\left| \sum_{n=1}^3 \frac{1}{2} g^{mn} 
\left(
\frac{\partial g_{ni}}{\partial x_j} + \frac{\partial g_{nj}}{\partial x_i} - \frac{\partial g_{ij}}{\partial x_n}
\right)
\right|^2\\&\leq
\sum_{i=1, j=1}^3 3^4 \delta^4 \\&\leq 
3^6 \delta^4. 
\end{split} 
\end{align}
In order to show the inequalities above we repeatedly use \cref{inverse metric bound} and our main assumption $\|{g-\gflat}_{C^1} \leq \delta$.
In the fourth step we use \cref{Christoffel Symbol Estimate}

\cref{nabla difference bound} implies that:
\begin{align*}
|\nabla^g dx_m -\nabla dx_m|_{\gflat} \leq 3 \Cr{Christoffel Bound} .
\end{align*}
By linearity we obtain:
\begin{align*}
|\nabla^g \omega- \nabla \omega|_{\gflat} \leq 3\sqrt{3}\Cr{Christoffel Bound}|\omega|_{\gflat}.
\end{align*}

Evidently, in order to prove \cref{Step 1: T3 Injectivity Estimate} we are interested in the above bound for the 1-form $\omega = du$
The bound we care about is the following:
\begin{align}
\label{nabla comparison estimate}
|\nabla^g du - \nabla du|_{\gflat}\leq 3 \Cr{Christoffel Bound} \sum_{i=1}^3 \left| \frac{\partial u}{\partial x_i} \right|,
\end{align}
combining \cref{nabla comparison estimate} and \cref{nabla estimate} gives us: 
\begin{align}
\label{Very Painful Estimate No More}
\begin{split}
|\nabla^g du|_{g} &\leq |\nabla^g du - \nabla du|_{g} + |\nabla du|_{g}
\\&\leq
\Cr{Co2Tensor Estimate}(|\nabla^g du - \nabla du|_{\gflat} + |\nabla du|_{\gflat}) \\&\leq 
\Cr{Co2Tensor Estimate}(3 \Cr{Christoffel Bound}\sum_{i=1}^3 \left| \frac{\partial u}{\partial x_i} \right| + |\nabla du|_{\gflat})
\\&\leq 
\Cr{Co2Tensor Estimate}(3\sqrt{3} \Cr{Christoffel Bound} \sqrt{\sum_{i=1}^{3} \left(\frac{\partial u}{\partial x_i} \right)^2} + |\nabla du|_{\gflat})
\\&\leq 
\Cr{Co2Tensor Estimate}(3 \sqrt{3} \Cr{Christoffel Bound} |du|_{\gflat} + |\nabla du|_{\gflat}).
\end{split}
\end{align}
In the penultimate step we use the equivalence of $1$-norm and Euclidean norm and in the last step we use the definition of $|du|_{\gflat}$. \\
We derive similarly: 
\begin{align}
\label{nabla df lower estimate}
\begin{split}
|\nabla^g du|_{g} &\geq \Cr{Co2Tensor Lower Estimate}(|\nabla du|_{\gflat}-  |\nabla^g du - \nabla du|_{\gflat})
\\&\geq
\Cr{Co2Tensor Lower Estimate} \left(|\nabla du|_{\gflat} - 3\Cr{Christoffel Bound} \sum_{i=1}^3 \left|\frac{\partial u}{\partial x_i}\right|\right) 
\\&\geq 
\Cr{Co2Tensor Lower Estimate} (|\nabla du|_{\gflat} - 3\sqrt{3} \Cr{Christoffel Bound} |du|_{\gflat}).
\end{split}
\end{align}

We now have gathered all the ingredients required to prove \cref{Step 1: T3 Injectivity Estimate}.
We start with \cref{Very Painful Estimate No More}, 
\begin{align*}
|\nabla^g du|_{g} \leq \Cr{Co2Tensor Estimate} (\sqrt{3} |du|_{\gflat} + |\nabla du|_{\gflat}).
\end{align*}
This pointwise bound readily gives a bound in the $L^p$ norm. 
\begin{align*}
\begin{split}
\|{\nabla^g du}_{g, L^p}
&= 
\left( \int_{T^3} |\nabla^g du|_{g} ^p \sqrt{det(g)} dx \right)^\frac{1}{p} 
\leq 
\Cr{Upper Determinant Estimate}^{\frac{1}{2p}} \Cr{Co2Tensor Estimate} 
\left\{
\left(
\int_{T^3} |\nabla du|_{\gflat}^p dx
\right)^\frac{1}{p} 
+ 
\left(
\int_{T^3} 3\sqrt{3} \Cr{Christoffel Bound} |du|_{\gflat}^p dx
\right)^\frac{1}{p}
\right\}. 
\end{split}
\end{align*}
Summing the above estimate with \cref{Step 1 Part 1: T3 Injectivity Estimate} gives us the result.
\begin{align*}
\begin{split}
\|{u}_{g, W^{2,p}(T^3)} 
&\leq
\Cr{S1P1 T3 Injectivity Estimate}[\left(\int_{T^3}|u|^p\right)^\frac{1}{p}+\left(\int_{T^3} |du|^p_{\gflat}\right)^\frac{1}{p}] + \Cr{Upper Determinant Estimate}^{\frac{1}{2p}} \Cr{Co2Tensor Estimate} \left(\int_{T^3} |\nabla du|^p_{\gflat} dx\right)^\frac{1}{p} + 3 \sqrt{3}\Cr{Upper Determinant Estimate}^{\frac{1}{2p}} \Cr{Christoffel Bound}\Cr{Co2Tensor Estimate}\left(\int_{T^3} |du|^p_{\gflat} dx\right)^\frac{1}{p}
\\&\leq
\left( \Cr{S1P1 T3 Injectivity Estimate}+ (1+3\sqrt{3}\Cr{Christoffel Bound}) \Cr{Co2Tensor Estimate} \Cr{Upper Determinant Estimate}^{\frac{1}{2p}} \right) \|{u}_{W^{2,p}}.
\end{split}
\end{align*}

In order to obtain the corresponding lower bound, we perform a similar calculation. We obtain the following integral estimate from \cref{nabla df lower estimate}:
\begin{align*}
\|{\nabla^g du}_{g, L^p} \geq \Cr{Co2Tensor Lower Estimate} \Cr{Lower Determinant Estimate}^{\frac{1}{2p}}  \left\{ \|{\nabla du}_{\gflat, L^p} - 3\sqrt{3} \Cr{Christoffel Bound} \|{du}_{\gflat, L^p}\right\},
\end{align*}
\begin{align*}
\|{u}_{g, W^{2,p}} \geq \Cr{S1P1 T3 Injectivity Estimate Lower Bound} (\|{u}_{\gflat, L^p} + \|{du}_{\gflat, L^p}) +  \Cr{Co2Tensor Lower Estimate} \Cr{Lower Determinant Estimate}^{\frac{1}{2p}} \left\{ \|{\nabla du}_{\gflat, L^p} - 3\sqrt{3} \Cr{Christoffel Bound} \|{du}_{\gflat, L^p}\right\}.
\end{align*}
\end{proof}

Our work in the first section immediately gives us the second step in \cref{overall strategy for injective estimate with arbitrary-ish riemannian metric}) 
This is the following:
\begin{align*}
  \|{u}_{W^{2,p}} \leq \Cr{const:injectivity-main-flat}\|{\Delta u}_{L^p}.
\end{align*}
In order to prove \cref{theorem:T3-injectivity-estimate}, it remains to compare Laplace operators of $\gflat$ and $g$.
\begin{lemma}
Assume $\|{g- \gflat}_{C^1} < \delta$, 
\begin{align}
\label{Delta Comparison Bound}
\left(\int_{T^3}  |\Delta u- \Delta^g u (\det g)^{\frac{1}{2p}}|^p\right)^{\frac{1}{p}} \leq \Cl{Absorption Laplacian Term} \|{u}_{W^{2,p}},
\end{align}
where $\Cr{Absorption Laplacian Term}$ is defined in \cref{Absorption Laplacian Term}.
\end{lemma}

\begin{proof}
We recall a convenient equivalent definition of $\Delta^g$ in terms of Christoffel Symbols. 
\begin{align*}
\Delta^g (u) = \sum_{i=1,j=1}^3 g^{ij} \left(\frac{\partial^2 u}{\partial x_i \partial x_j} - \sum_{k=1}^3 \Gamma^k_{ij} \frac{\partial u}{\partial x_k}\right).
\end{align*}
We now compute that:
\begin{align*}
\begin{split}
\|{\Delta u - \Delta^g u (\det g)^{\frac{1}{2p}}}_{\gflat ,L^p} 
&\leq \|{\sum_{i=1,j=1}^3 (g^{ij} (\det g)^{\frac{1}{2p}} - \delta_{ij}) \frac{\partial^2 u}{\partial x_i \partial x_j}}_{L^p} + \|{\sum_{i=1. j=1, k=1}^3 g^{ij}\Gamma^k_{ij} (\det g)^{\frac{1}{2p}}\frac{\partial u}{\partial x_k}.}_{L^p}
\\&\leq 
\left(\left(\Cr{Upper Determinant Estimate}^{\frac{1}{2p}}-1\right)(1+2\delta) + 2\delta\right) \|{\sum_{i=1, j=1}^3 \frac{\partial^2 u}{\partial x_i x_j}}_{\gflat, L^p} + 9 (1+2\delta) \Cr{Christoffel Bound} \Cr{Upper Determinant Estimate}^{\frac{1}{2p}} \|{\sum_{k=1}^3 \frac{\partial u}{\partial x_k}}_{L^p} 
\\&\leq
3\left(\left(\Cr{Upper Determinant Estimate}^{\frac{1}{2p}}-1\right)(1+2\delta) + 2\delta\right) \|{D^2 u }_{\gflat, L^p} + 9 \sqrt{3} (1+2 \delta) \Cr{Christoffel Bound} \Cr{Upper Determinant Estimate}^{\frac{1}{2p}} \|{Du}_{L^p} 
\end{split}
\\&\leq
3\left(\left(\Cr{Upper Determinant Estimate}^{\frac{1}{2p}}-1\right)(1+2\delta) + 2\delta\right)+9 \sqrt{3} (1+ 2\delta)\Cr{Christoffel Bound} \Cr{Upper Determinant Estimate}^{\frac{1}{2p}}) \|{u}_{W^{2,p}}.
\end{align*}
In the second step we use the bound presented below on the summand exclusively in i and j, for the latter term involving Christoffel Symbols we use \cref{Christoffel Symbol Estimate} and in the third step we use norm equivalence of the 1-norm and the 2-norm. 
\begin{align*}
|g^{ij} (\det g)^{\frac{1}{2p}} - \delta_{ij}|\leq |g^{ij} (\det g)^{\frac{1}{2p}} - g^{ij}| + |g^{ij}- \delta_{ij}| \leq \left(\Cr{Upper Determinant Estimate}^{\frac{1}{2p}}-1\right)(1+2\delta) + 2\delta.
\end{align*}
The factor of 2 in the first summand in the second step arises due  to \cref{inverse metric bound}.
We conclude that 
\begin{align}
\label{Absorption Laplacian Term}
\Cr{Absorption Laplacian Term}(\delta)=\Cr{Absorption Laplacian Term} = 3\left(\left(\Cr{Upper Determinant Estimate}^{\frac{1}{2p}}-1\right)(1+2\delta) + 2\delta\right) + 9 (1+2\delta)\Cr{Christoffel Bound} \sqrt{3} \Cr{Upper Determinant Estimate}^{\frac{1}{2p}}.
\end{align}
\end{proof}

We are now ready to prove the main result of the section.

\begin{proof}[Proof of \cref{theorem:T3-injectivity-estimate}]
First observe that
\begin{align}
\label{Triangle Inequality Delta}
\|{\Delta u}_{\gflat, L^p} \leq \|{\Delta u - \Delta^g u (\det g)^{\frac{1}{2p}}}_{\gflat, L^p} + \|{\Delta^g u}_{g, L^p}.
\end{align}
 We consider \cref{proposition:injectivity-estimate-flat-metric} ,substituting \cref{Triangle Inequality Delta} and \cref{Delta Comparison Bound} into the right hand side, we obtain:
\begin{align*}
\|{u}_{\gflat, L^p_2} \leq \Cr{const:injectivity-main-flat} \|{\Delta^g u}_{g, L^p} + \Cr{const:injectivity-main-flat} \Cr{Absorption Laplacian Term} \|{u}_{\gflat, L^p_2},
\end{align*}
by taking $\delta >0$ to be sufficiently small, $\Cr{Absorption Laplacian Term} < \frac{1}{\Cr{const:injectivity-main-flat}}$, it follows:
\begin{align*}
\|{u}_{\gflat, L^p_2} \leq \frac{\Cr{const:injectivity-main-flat}}{1-\Cr{const:injectivity-main-flat}\Cr{Absorption Laplacian Term}}\|{\Delta^g u}_{g, L^p}
\end{align*}
Using \cref{Step 1: T3 Injectivity Estimate} we obtain the required result with:
\begin{align}
\label{equation:injectivity-nonflat-constant-def}
\Cr{const:injectivity-main-nonflat}= \frac{\Cr{Step 1 3.3}\Cr{const:injectivity-main-flat}}{1- \Cr{const:injectivity-main-flat}\Cr{Absorption Laplacian Term}}.
\end{align}
\end{proof}

\begin{remark}
Using Wolfram Mathematica we compute that in fact $\Cr{const:metric-max-difference} < 3 \times 10^{-14} $ satisfies the condition from \cref{theorem:T3-injectivity-estimate}.
\end{remark}

\section{Application: nowhere vanishing harmonic 1-forms}
\label{section:application-nowhere-vanishing-1-forms}


We are now ready to give the proof of \cref{theorem:application-nowhere-vanishing-1-form}:

\begin{proof}[Proof of \cref{theorem:application-nowhere-vanishing-1-form}]
    Let $\xi := \Delta_g^{-1} (\d_g ^* \d \,x_1)$, then
\begin{align}
\begin{split}
\Delta_g(d x_1 - d\xi) &= \Delta^g \{d x_1- d\Delta_g^{-1}(d^*_g dx_1)\}
                 \\&= \Delta^g d x_1 - \Delta^g d \Delta_g^{-1} (d^*_g d x_1)
                 \\&= \Delta^g d x_1 - dd^*_g d x_1
                 \\&= \Delta^g d x_1 - dd^*_g \d x_1 - d^*_g d d x_1=0. 
\end{split}
\end{align}
In the second step we use the fact that d and $\Delta^g$ commute and in the third step we use the fact that $\d^2=0$.
We now verify that the one-form $\d x_1 - \d\xi$ is nowhere-vanishing. 

    The function $\xi$ satisfies
    \begin{align*}
        \|{\d \xi}_{g, C^0}
        &\leq
        \Cr{Covector Estimate}\|{\d \xi}_{C^0,{\gflat}}
        \\&\leq
        \Cr{Morrey} \Cr{Covector Estimate}
        \|{\xi}_{\gflat, W^{2,4}}
        \\&\leq
        \frac{\Cr{Morrey} \Cr{Covector Estimate}}{\Cr{Step 1 3.3 Lower Bound}}
        \|{\xi}_{g, L^4_2}
        \\&\leq
        \frac{\Cr{Morrey} \Cr{Covector Estimate}\Cr{const:injectivity-main-nonflat}}{\Cr{Step 1 3.3 Lower Bound}}\|{\d_g ^* \d x_1}_{g, L^4},
    \end{align*}
    where we use \cref{comparison lemma covectors} in the first step,
    we use \cref{corollary:morrey's-ineq} in the second step, 
    we use \cref{theorem:L^p_2-norm-comparison} in the third step,
    and use \cref{theorem:T3-injectivity-estimate} in the last step.
    \\
    Denote:
    \begin{align}
        \Cl{1-Form Injectivity Estimate}= \frac{\Cr{Morrey} \Cr{Covector Estimate}\Cr{const:injectivity-main-nonflat}}{\Cr{Step 1 3.3 Lower Bound}}.
   \end{align}

We use a definition of the codifferential which involves Christoffel Symbols, for convenience of calculation.
For $\omega=\sum_j \omega_j \d x_j$:
\begin{align}
d^*_g \omega=-\sum_{i,j}g^{ij} \partial_i \omega_j + \sum_{i,j,k} g^{ij} \Gamma^k_{ij} \omega_k,
\end{align}
letting $\omega = \d x_1$ in the above expression gives:
\begin{align}
d^*_g \d x_1 = \sum_{i,j} g^{ij} \Gamma^1_{ij},
\end{align}
because each of the $\omega_j$ are constant functions.
\cref{lemma:christoffel-bound,inverse metric bound} then give:
\begin{align}
|d^*_g \d x_1| \leq (3+18\delta)\Cr{Christoffel Bound}.
\end{align}

    This then implies the claim, because for all $x \in T^3$ we have that
    \begin{align}
    \begin{split}
        |\d \, x_1+\d \, \xi|_{g}(x)
        &\geq
        |\d \, x_1|_{g}(x)
        -
        |\d \, \xi|_{g}(x)
        \\&\geq
        \Cr{Covector Lower Estimate}-\|{\d \, \xi}_{g,C^0}
        \\&\geq
        \Cr{Covector Lower Estimate}- \Cr{1-Form Injectivity Estimate} \|{\d^*_g \d x_1}_{g, L^p}
        \\&\geq
        \Cr{Covector Lower Estimate}- \Cr{1-Form Injectivity Estimate} (3+18\delta)\Cr{Christoffel Bound} \Cr{Upper Determinant Estimate}^{\frac{1}{2p}},
    \end{split}    
    \end{align}
    which is $>0$ for small $\delta$.
    That is, the $1$-form $\d x_1+\d \, \xi \in \Omega^1(T^3)$ is nowhere vanishing.
    For future use, we denote:
    \begin{align}
    \label{End of Proof}
    \Cl{endofpaper}= \Cr{endofpaper}(\delta)=   \Cr{Covector Lower Estimate}- \Cr{1-Form Injectivity Estimate} (3+18\delta)\Cr{Christoffel Bound} \Cr{Upper Determinant Estimate}^{\frac{1}{2p}}.
    \end{align}
\end{proof}

\appendix
\section{Appendix}

\subsection{Marcinkiewicz interpolation theorem}

\begin{theorem}[Theorem 9.8 in \cite{Gilbarg1998}]
    \label{theorem:marcinkiewicz}
    Let $T$ be a linear mapping from $L^1(\Omega) \cap L^r(\Omega)$ into itself, $1 \leq q < r < \infty$, and suppose there are constants $T_1$ and $T_2$ such that
    \begin{align}
        \mu_{Tf}(t)
        \leq
        \left(
        \frac{T_1 \|{f}_{L^q}}{t}
        \right)^q,
        \quad
        \mu_{Tf}(t)
        \leq
        \left(
        \frac{T_2 \|{f}_{L^r}}{t}
        \right)^r,
    \end{align}
    for all $f \in L^q(\Omega) \cap L^r(\Omega)$ and $t>0$.
    Then $T$ extends as a bounded linear mapping from $L^p(\Omega)$ into itself for any $p$ such that $q < p<r$, and
    \begin{align}
        \|{Tf}_{L^p}
        \leq
        C_{\text{Marcinkiewicz}}(p,q,r) T_1^\alpha T_2^{1-\alpha} \|{f}_{L^p},
    \end{align}
    where $\frac{1}{p} = \frac{\alpha}{q}+\frac{1-\alpha}{r}$ and
    \[
    C_{\text{Marcinkiewicz}}(p,q,r)
    =
    2
    \left(
    \frac{p(r-q)}{(p-q)(r-p)}
    \right)^{1/p}.
    \]
\end{theorem}

\begin{proof}
The same proof can be found in \cite{Gilbarg1998}, to make this paper more self-contained we present the proof here also.
These are two basic facts about distribution functions which are also defined below

For any $p>0$ and $|u|^p \in L^1(\Omega)$:
\begin{align}
\label{Markov Inequality}
\mu_{u}(t):= |\{x\in \Omega | u(x) > t\}| \leq \frac{\int_{\Omega} |f|^p}{t^p}
\end{align}
\begin{align}
\label{Fubini Lower Moment Estimate}
\int_{\Omega}|u|^p = p \int_{0}^{\infty} t^{p-1} \mu_{u}(t) dt.
\end{align}

Fix $u \in L^q(\Omega) \cap L^r(\Omega)$ and s>0, then rewrite $u = f_1 +f_2$
where
\begin{align}
f_1(x):=\begin{cases}
         u(x) & \text{if }  |u(x)|>s\\
         0 & \text{if }  |u(x)|\leq s
        \end{cases}
\\
f_2(x):=\begin{cases}
         0 & \text{if }  |u(x)|>s\\
         u(x) & \text{if }  |u(x)|\leq s
        \end{cases}.   
\end{align}
Using the triangle inequality and the assumptions we obtain that:
\begin{align}
\mu_{Tu}(t) &\leq \mu_{Tf_1}\left(\frac{t}{2}\right) +\mu_{Tf_2}\left(\frac{t}{2}\right)\\
            &\leq \left(\frac{2T_1}{t}\right)^q \int_{\Omega} |f_1|^q + \left(\frac{2T_2}{t}\right)^r \int_{\Omega} |f_2|^r.
\end{align}
We combine the above with \cref{Fubini Lower Moment Estimate} to obtain:
\begin{align}
\int_{\Omega}|Tu|^p &= p\int_{0}^{\infty} t^{p-1}\mu_{Tf}(t) dt\\
                    &= p(2T_1)^q\int_{0}^{\infty} t^{p-q-1}\left(\int_{\Omega} |f_1|^q\right) dt  + p(2T_2)^r \int_{0}^{\infty} t^{p-r-1}\left(\int_{\Omega} |f_2|^r\right)dt\\
                    &= p(2T_1)^q\int_{0}^{\infty} t^{p-q-1}\left(\int_{|f|>s} |f|^q\right) dt  + p(2T_2)^r \int_{0}^{\infty} t^{p-r-1}\left(\int_{|f|\leq s} |f|^r\right)dt.
\end{align}
The main trick presented in the proof \cite[Theorem 9.8]{Gilbarg1998} is that we take $t=As$ where A is a positive constant to be fixed later. We first derive our required estimate with the constant being a function of A and then minimise it to obtain a sharper result. 
\begin{align}
\int_{\Omega}|Tu|^p\leq p(2T_1)^q A^{p-q}\int_{0}^{\infty}s^{p-q-1}\left(\int_{|f|> s} |f|^q\right) ds + p(2T_1)^r A^{p-r}\int_{0}^{\infty} s^{p-r-1}\left(\int_{|f|\leq s} |f|^r\right) ds.
\end{align}
The above integrals on the RHS have strictly non-negative integrand and thus can be rewritten and evaluated using Fubini's Theorem as follows:
\begin{align}
  \int_{0}^{\infty} s^{p-q-1}\left(\int_{|f| > s} |f|^q\right) ds 
   &= \int_{\Omega} |f|^q \int_{0}^{|f|} s^{p-q-1} ds \\
   &= \frac{1}{p-q} \int_{\Omega} |f|^p,
\end{align}
\begin{align}
   \int_{0}^{\infty} s^{p-r-1}\left(\int_{|f|\leq s} |f|^r\right) ds
   &= \int_{\Omega} |f|^r \int_{|f|}^{\infty} s^{p-r-1} ds\\
   &= \frac{1}{r-p}\int_{\Omega} |f|^p.
\end{align}

Substituting these evaluated integrals above gives us the required bound as a function of $A$ as promised. 
\begin{align}
\int_{\Omega}|Tu|^p \leq \{\frac{p(2T_1)^q A^{p-q}}{p-q} + \frac{p(2T_1)^r A^{p-r}}{r-p}\} \int_{\Omega} |u|^p.
\end{align}
Minimising the expression in braces with respect to A gives us the required result. 
\end{proof} 
\subsection{Explicit Sobolev Constant}
\begin{theorem}[Explicit Sobolev Inequality for compactly supported functions]
    \label{Sobolev Embedding}
         If $1\leq q<n$, then all $u \in W^{1, q}_0( \R ^n) \text{ satisfy}:$
    \begin{equation}
         \|{u}_{L^p} \leq K(n,q)\|{Du}_{L^q}, 
    \end{equation}
with $\frac{1}{p} = \frac{1}{q} -\frac{1}{n}$ and
\begin{align*}
K(n,q)= \frac{q-1}{n-q}\left(\frac{n-q}{n(q-1)}\right)^\frac{1}{q} \left(\frac{\Gamma(n+1)}{\Gamma(n/q)\Gamma(n+1-\frac{n}{q})\omega_{n-1}}\right)^\frac{1}{n},
\end{align*}
where $\omega_{n}$ is the volume of the n-dimensional sphere with radius 1. \cite[2.14]{Aubin1998}
\end{theorem}
\begin{remark}
We use the above claim in \cref{Middle Term Estimate} with $n=3, q=2$, denote $\Cl{sobolev-embedding} = K(3,2)$
\end{remark}
\begin{theorem}[Explicit Sobolev Inequality (full generality) on the cube]
\label{other-sobolev-embedding}
Suppose $p=4$, $\Omega = [0,1]^3$ and $u \in W^{1,2}(\Omega)$ , Then the following inequality holds:
\begin{align}
\|{u}_{L^4(\Omega)} \leq \Cl{Sobolev-Embedding-Cube} \|{u}_{W^{1,2}(\Omega)}.
\end{align}
\end{theorem}
This can be found in \cite[page 15, table 6]{MizuguchiTanakaSekineOishi2017} with $\Cr{Sobolev-Embedding-Cube} = 13.25$ 
\subsection{Poincaré Inequality}
We first prove the following lemma.
\begin{theorem}[Poincaré Inequality]
\label{Poincaré Inequality}
Suppose $u\in W^{1,p}_0(\Omega)$ and $1\leq p < \infty$, then we have the following inequality:
\begin{align}
\|{u}_{L^p}
\leq 
\left(
\frac{1}{\omega_n} |\Omega|
\right)^\frac{1}{n} \|{Du}_{L^p}.
\end{align}
\end{theorem}

There are two tools which we use to prove this theorem, one represents $u(x)$ as an integral (specifically a Riesz Potential to be defined below) and then refer to a bound for the norm of Riesz Potential:

\begin{lemma}[Lemma 7.14 \cite{Gilbarg1998}]
  \label{Non-trivial Change Of Variables}
  Let $u\in W^{1,1}_0 (\Omega)$ Then:
  \begin{align}
  u(x) = \frac{1}{n\omega_n} \sum_{i=1}^n\int_{\Omega} \frac{(x_i-y_i)D_iu(y)}{|x-y|^n}dy. 
  \end{align} 
\end{lemma}
\begin{proof}
We use the same proof as \cite[Lemma 7.14]{Gilbarg1998} but we provide the elementary details which are omitted.
First we prove the claim for $u \in C^1_1(\Omega)$:
For an arbitrary choice of direction $\omega \in \R^n$ and $|\omega|=1$ we have that:
\begin{align}
u(x) &= - \int_{0}^{\infty} \frac{d}{dr}(u(x+r\omega))dr 
   \\&= -\int_0^\infty D_{\omega}(u(x+r\omega)) dr
   \\&= -\frac{1}{n\omega_n}\int_{0}^{\infty} \int_{|\omega|=1} D_{\omega}(u(x+r\omega))dr
   \\&=  -\frac{1}{n\omega_n}\int_{0}^{\infty} \int_{|\omega|=1} \omega \cdot D (u(x+rw))dr
   \\&= -\frac{1}{n \omega_n}\int_{\Omega} \frac{(y-x) \cdot D(u(y))}{|y-x|} \frac{1}{|y-x|^{n-1}} dy
   \\&= \frac{1}{n\omega_n} \sum_{i=1}^n\int_{\Omega} \frac{(x_i-y_i)D_iu(y)}{|x-y|^n}dy
\end{align}
For the first inequality we use the fundamental theorem of calculus (and the fact that $u$ is compactly supported), for the second we use the multivariate chain rule (noting that $D_{\omega}(u(x+r\omega)) = \omega \cdot D(u(x+r\omega))$, for the third we integrate over $|\omega|=1$ and divide by the measure of this set (surface area of the n-1 sphere) (one integrates LHS and RHS to get this), for the fifth equality we change from spherical coordinates back to cartesian.

\end{proof}
\begin{remark}
By applying the Cauchy-Schwarz Inequality to the above claim we obtain a useful estimate. 
\begin{align}
\label{compact support representation}
|u(x)|\leq \frac{1}{n\omega_n}\int_{\Omega} \frac{|Du(y)|}{|x-y|^{n-1}} dy.
\end{align}
\end{remark}
\begin{lemma}[Bound for norm of Riesz Potential, Lemma 7.12 in \cite{Gilbarg1998}]
    \label{Riesz Potential Bound}
    We define the operator $V_{\mu}$ on $L^1(\Omega)$ to be the Riesz Potential:
    \begin{align}
    (V_{\mu}f)(x) = \int_{\Omega} |x-y|^{n(\mu-1)} f(y) dy
    \end{align}
    The operator $V_{\mu}$ maps $L^p(\Omega)$ continuously into $L^q(\Omega)$ for any $1\leq q < \infty$ satisfying 
    \begin{align}
    0\leq \delta = \delta(p,q) = p^{-1} -q^{-1} < \mu 
    \end{align}
     Additionally, for any $f\in L^p(\Omega)$ :
    \begin{align}
    \|{V_\mu f}_{L^q} \leq 
    \left(
    \frac{1-\delta}{\mu-\delta}
    \right)
    ^{1-\delta} \omega_n^{1-\mu} |\Omega|^{\mu - \delta}\|{f}_{L^p}.
    \end{align}
\end{lemma}

\begin{proof}[Proof of \cref{Poincaré Inequality}]
We are now ready to prove the Poincaré Inequality.
\begin{align}
  |u(x)|&= |\frac{1}{n \omega_n}\sum_{i=1}^{n}\int_{\Omega} \frac{(x_i-y_i)D_iu(y)}{|x-y|^n} dy|
  \\&\leq \frac{1}{n\omega_n} \int_{\Omega} \sum_{i=1}^n \frac{|x_i -y_i| |D_i u(y)|}{|x-y|^{n}} dy 
  \\&\leq \frac{1}{n\omega_n} \int_{\Omega} \frac{(\sum_{i=1}^{n} |x_i - y_i|^2)^\frac{1}{2}(\sum_{i=1}^{n} |D_iu(y)|^2)^\frac{1}{2}}{|x-y|^n} dy 
  \\&= \frac{1}{n \omega_n}\int_{\Omega} \frac{|Du(y)|}{|x-y|^{n-1}} dy
  \\&= (V_{\frac{1}{n}} |Du|)(x).
\end{align}
For the first equality we use \cref{Non-trivial Change Of Variables}, for the second inequality we use the triangle inequality, for the third we use the Cauchy-Schwarz inequality, the last equality is using the definition of Riesz potential.
\begin{align}
\|{u}_{L^p} \leq \frac{1}{n \omega_n} n \omega_n^{1-\frac{1}{n}} |\Omega|^\frac{1}{n} \|{Du}_{L^p}= \left(\frac{1}{\omega_n} |\Omega|\right)^\frac{1}{n} \|{Du}_{L^p}. 
\end{align}
We take the $L^p$ norm of the derived inequality above and then we apply \cref{Riesz Potential Bound}, this concludes the proof.
\end{proof}

\subsection{Morrey's Inequality}
In this section we aim to prove the following:
\begin{theorem}
\label{Morrey's Inequality}
Suppose $u \in C^1(T^3) $, then for $p=4 $, we have the following inequality.
\begin{align}
\|{u}_{C^0(T^3)} \leq \Cl{Morrey} \|{u}_{W^{1,p}(T^3)}.
\end{align}
\end{theorem}
\begin{proof}
It is sufficient to prove the same claim for periodic functions defined on $\Omega = [0,1]^3$.
A more general claim is proved in \cite[Theorem 4, Section 5.6]{evans2010partial}, the proof is constructive and we obtain the desired constant following this construction.  
The rough idea is as follows:
A similar bound to \cref{compact support representation} can be exhibited for mean zero functions rather than compactly supported functions, this integral representation readily gives the result. This heuristic is made precise below.

\begin{align}
\label{Morrey Setup}
\begin{split}
|u(x)| &= \frac{1}{\omega_3} \int_{B(x,1)} |u(x)| dy \\&
       \leq \frac{1}{\omega_3} \int_{B(x,1)} |u(x)-u(y)| dy + \frac{1}{\omega_3}  \int_{B(x,1)} |u(y)| dy.
\end{split}
\end{align}
We first prove the following:
\begin{lemma}
\label{Morrey Spheres}
\begin{align}
\int_{B(x,r)} |u(x)-u(y)| dy \leq \frac{r^n}{n} \int_{B(x,r)} \frac{|Du (y)|}{|x-y|^{n-1}} dy. 
\end{align}
\end{lemma}
The proof of this lemma is \cite[Theorem 4, Section 5.6]{evans2010partial}. We follow the proof carefully and obtain explicit constants for the bounds.
We aim to apply Hölder's inequality to the right hand side of \cref{Morrey Spheres} in order to obtain a $\|{Du}_{L^p}$ term. Applying Hölder's inequality to the second integral in \cref{Morrey Setup} will yield the required claim with explicit constants for spheres.

We first fix $0<s<r$ and $\omega \in \mathbb{S}^n$, we obtain that:
\begin{align}
\begin{split}
|u(x+s\omega) - u(x)| 
&= 
\left| \int_{0}^s \frac{d}{dt} u(x+t\omega) dt \right|
= 
\left| \int_{0}^s Du(x+t\omega)\cdot \omega dt
\right|
\leq 
\int_{0}^s |Du(x+t\omega)| dt .
\end{split}
\end{align}
In the first step we use the fundamental theorem of calculus, in the second we use the multivariate chain rule and in the third we use that $\omega$ has length one. 
\\
Using Fubini's Theorem we obtain that:
\begin{align}
\int_{\omega \in \mathbb{S}^n} |u(x+s \omega) - u(x)|d\omega \leq \int_{0}^s \int_{\omega \in \mathbb{S}^n} |Du(x+t\omega)| d\omega dt.
\end{align}
We now use the substitution $y = x+ t \omega$ to convert the above integral from spherical to cartesian coordinates.
\begin{align}
\begin{split}
\int_{\omega \in \mathbb{S}^n} |u(x+s\omega) - u(x)| d\omega 
&\leq \int_{y \in B(x,s)} \frac{|Du(y)|}{|x-y|^{n-1}} dy 
\\&\leq \int_{y \in B(x,r)} \frac{|Du(y)|}{|x-y|^{n-1}} dy .
\end{split}
\end{align}
In the integral on the left hand side we use the substitution $z= x+ s\omega$, we obtain the following:
\begin{align}
\int_{\omega \in \mathbb{S}^n} |u(x+s\omega) - u(x)| d\omega = \frac{1}{s^{n-1}}\int_{z \in \partial B(x,s)} |u(z) - u(x)| dz,  .
\end{align}
thus we have that:
\begin{align}
\int_{z \in \partial B(x,s)} |u(z) -u(x)| dz \leq s^{n-1} \int_{y \in B(x,r)} \frac{|Du(y)|}{|x-y|^{n-1}} dy. 
\end{align}
Integrating the left and right hand side from 0 to r with respect to s yields the claim.
\begin{align}
\int_{B(x,r)} |u(x)-u(y)| dy \leq \frac{r^n}{n} \int_{B(x,r)} \frac{|Du (y)|}{|x-y|^{n-1}} dy
\end{align}
Substituting \cref{Morrey Spheres} into \cref{Morrey Setup} with $r=1$ yields:
\begin{align}
|u(x)| \leq \frac{1}{\omega_3} \int_{B(x,1)}|u(y)| dy + \frac{1}{n \omega_3}\int_{B(x,1)} \frac{|Du(y)|}{|x-y|^{n-1}} dy,
\end{align}
Applying Hölder's Inequality to both integrals yields, using $p=4$, the following:
\begin{align}
|u(x)| \leq \frac{1}{\omega_3}  (\omega_3)^\frac{3}{4} \|{u}_{L^p(B(x,1))} + 
\frac{1}{n \omega_3}
\left(
\int_{B(x,1)}|x-y|^{\frac{4}{3}(n-1)}
\right)^{\frac{3}{4}} \|{Du}_{L^p(B(x,1))}.
\end{align} 
The unit ball $B(x,1)$ contains a translate of $\Omega=[0,1]^3$. 
By periodicity we arrive at the following:
\begin{align}
\begin{split}
|u(x)| &\leq \omega_3^{-\frac{1}{4}} \|{u}_{L^p(T^3)} + 
\left(
\int_{B(x,1)}|x-y|^{\frac{4}{3}(n-1)}
\right)^{\frac{3}{4}} \|{Du}_{L^p(T^3)}
\\&\leq
\max{(\omega_3^{\frac{1}{4}}, (\int_{B(x,1)}|x-y|^{\frac{4}{3}(n-1)})^{\frac{3}{4}}} \|{u}_{L^p_1(T^3)}.
\end{split}
\end{align}
We compute that in the case that n=3: $\Cr{Morrey}=\max{(\omega_3^{-\frac{1}{4}}, (\int_{B(x,1)}|x-y|^{\frac{8}{3}})^{\frac{3}{4}}} \|{u}_{L^p_1(T^3)} $
\end{proof}

We prove a corollary that is used in \cref{section:application-nowhere-vanishing-1-forms}
\begin{corollary}
\label{corollary:morrey's-ineq}
Suppose $u \in L^p_2(T^3)$, then the following inequality holds:
\begin{align}
\|{\d u }_{C^0} \leq \Cr{Morrey}\|{\d u}_{L^p_1}. 
\end{align}
\end{corollary}
The claim follows upon computation of the pointwise norm $|\d u |_{g_{\text{flat}}}$
\begin{proof}
\begin{align}
\begin{split}
\|{\d u}_{\gflat, C^0} &\leq \|{\sqrt{\sum_{i=1}^3 \left(\frac{\partial u}{\partial x_i}\right)^2}}_{\gflat, C^0} \\&\leq
\Cr{Morrey} \|{\sqrt{\sum_{i=1}^3 \left(\frac{\partial u}{\partial x_i}\right)^2}}_{g_{\text{flat}, L^p_1}} \\&=
\Cr{Morrey} \|{\d u}_{g_{\text{flat}}, L^p_1}. 
\end{split}
\end{align}
\end{proof}
\section{Cut-off functions}

The function $S$ defined below was introduced in \cite{Perlin2002} as the \emph{smootherstep} function.
It interpolates between $0$ and $1$ and is twice continuously differentiable.

\begin{proposition}
    \label{proposition:cut-off-estimates}
    Let
    \begin{align}
        \begin{split}
            S: \R & \rightarrow [0,1]
            \\
            t & \mapsto 
            \begin{cases}
                0 &\text{ if } t <0
                \\
                6t^5-15t^4+10t^3 & \text{ if } t \in [0,1]
                \\
                1 & \text{ if } t>1
            \end{cases}
        \end{split}
    \end{align}
    define $Q=[0,1]^3 \subset \R^3$ and $\tilde{Q}=[-1,2]^3 \subset \R^3$ and let
    \begin{align}
        \begin{split}
            \chi: \tilde{Q} & \rightarrow [0,1]
            \\
            x & \mapsto d(x,Q)=\min_{q \in Q} |x-q|.
        \end{split}
    \end{align}
    Then $\chi$ is twice continuously differentiable and
    \begin{align}
        \begin{split}
            \|{\frac{\partial^2}{\partial x^2}\chi}_{L^\infty(\tilde{Q})}
            &\leq
            \Cl{partial-ii-chi-bound},
            \\
            \|{\frac{\partial}{\partial x} \chi}_{L^\infty(\tilde{Q})}
            &\leq
            \Cl{partial-i-chi-bound}
            \\
            \|{\frac{\partial^2}{\partial x\partial y} \chi}_{L^\infty(\tilde{Q})}
            &\leq
            \Cl{partial-ij-chi-bound}.
        \end{split}
    \end{align}
    for
    \[
        \Cr{partial-ii-chi-bound}
        =
        -60(-6+ \frac{13}{6}(9+ \sqrt{3}) -\frac{1}{4}(9+\sqrt{3})^2 +\frac{1}{108}(9+\sqrt{3})^3),
        \quad
        \Cr{partial-i-chi-bound}
        =
        60(2\sqrt{3}-3),
        \quad
        \Cr{partial-ij-chi-bound}
        =
        20(5 \sqrt{3}-6 ).
    \]
The above inequalities hold up to permutation of x,y and z by the symmetry of the cut off function. 
We use these bounds to compute $\Cr{delta-chi-bound}$ and $\Cr{D-chi-bound}$
\begin{align}
\begin{split}
\|{\Delta \chi}_{L^\infty(\tilde{Q})} &\leq 3\Cr{partial-ii-chi-bound},
\\
\|{D\chi}_{L^\infty(\tilde{Q})} &\leq \sqrt{3} \Cr{partial-i-chi-bound},
\\
\|{D^2\chi}_{L^\infty(\tilde{Q})} &\leq 3 \Cr{partial-ij-chi-bound}.
\end{split}
\end{align}
We conclude that
\begin{align}
\begin{split}
\Cr{delta-chi-bound} &= 3 \Cr{partial-ii-chi-bound},
\\
\Cr{D-chi-bound} &= \sqrt{3} \Cr{partial-i-chi-bound},
\\
\Cr{D^2-chi-bound} &= 3 \Cr{partial-ij-chi-bound}.
\end{split}
\end{align}
\end{proposition}

\section{Index of Constants}
In this section we list all the constants used in this paper and where they arise.
When the value of a certain constant is explicitly stated, it is larger than the constant computed by Mathematica. When computing future constants in terms of prior constants we use the values stored in Mathematica and then round up. 
\begin{subsection}{\cref{subsection:local-estimates}, \cref{subsection:estimates-on-t3-with-flat-metric}}
\begin{itemize}
\item $\Cr{const:injectivity-main-flat} = \frac{1}{4\pi^2}\{\Cr{constant:Schauder-Estimate} \Cr{Sobolev-Embedding-Cube}(1 + 27\sqrt{17} )\} + \Cr{constant:Schauder-Estimate}(1+27\Cr{Sobolev-Embedding-Cube})  $, The main result of Section 3.1: injectivity estimate with respect to the flat metric, see \cref{proposition:injectivity-estimate-flat-metric}
\item $\Cr{constant:mu-t-L2-norm}, \Cr{constant:mu-t-L1-norm}$, bounds on Distribution Functions which allow us to prove the Calderon-Zygmund Estimate. Applying \cref{theorem:marcinkiewicz} to \cref{Distribution Function Bounds} proves  \cref{theorem:calderon-zygmund-potential-formulation}
\item $\Cr{constant:mu-tg-L1-norm}, \Cr{constant:mu-tb-L1-norm}$, these bounds are used to derive $\Cr{constant:mu-t-L1-norm}$
\item $\Cr{constant:-2nd-Derivative-Bound-Newtonian-Potential} $ see \cref{lemma:2nd-Derivative-Newtonian} and \cite[Lemma 4.4]{Gilbarg1998}
\item $\Cr{constant:Integral-Bound-Calderon-Zygmund}$, used in the proof of \cref{theorem:calderon-zygmund-potential-formulation} 
\item $\Cr{constant:-F^*-measure-bound}$, in the proof of \cref{theorem:calderon-zygmund-potential-formulation}, see \cref{balls and cubes} and \cref{dist F* bound}
\item\[C_{Calderon-Zygmund}(n,p) = \begin{cases}
            C_{Marcinkiewicz}(1,2,p')T_1^\alpha T_2^{1-\alpha} = 2(\frac{p'-2}{1-p'})^{1/p}T_1^\alpha T_2^{1-\alpha}& p>2
            \\
            C_{Marcinkiewicz}(1,2,p) T_1^\alpha T_2^{1- \alpha}= 2(\frac{p-2}{1-p})^{1/p}T_1^\alpha T_2^{1-\alpha} & p<2
            \\
            1 & p=2
            \end{cases}\]
see \cref{theorem:calderon-zygmund-potential-formulation}
\item $\Cr{constant:Calderon-Zygmund-3.1} = C_{Calderon-Zygmund}(3,4) = 293.519$, this particular case of the Calderon-Zygmund estimate is used throughout the rest of \cref{subsection:estimates-on-t3-with-flat-metric} and \cref{subsection:estimates-on-t3-with-non-flat-metric},  see \cref{corollary:Calderon-Zygmund-3.1} 
\item $\Cr{D+D^2-regularity-estimate}$, see \cref{corollary:D+D^2-regularity-estimate}, this shows $L^2$ orthogonality to $Ker \Delta$ in \cref{theorem:T3-injectivity-estimate} is only necessary to bound the norm of $\|{f}$
\item $\Cr{Poincaré Inequality} = (\frac{1}{\omega_n} |\Omega|)^\frac{1}{n}$, an explicit estimate for the Poincaré Inequality, see \cref{Poincaré Inequality}, is dependent on $\Omega \subset \R^n$
\item $\Cr{Poincaré-Inequality-Q-tilde} = \left(\frac{81}{4\pi}\right)^{\frac{1}{3}}$, the explicit value for the Poincaré Inequality \cref{Poincaré Inequality} when $\Omega = \tilde{Q}$. 
\item $\Cr{constant:Schauder-Estimate} =  27(1+\Cr{D+D^2-regularity-estimate}\Cr{delta-chi-bound} + (2+54(27)\sqrt{17}) \Cr{D+D^2-regularity-estimate}\Cr{Holder Q-Qtilde} \Cr{D^2-chi-bound} \Cr{sobolev-embedding}$, denotes an explicit value for which the Schauder Estimate holds on $T^3$, see \cref{equation:schauder-estimate-on-T3} 
\item $\Cr{delta-chi-bound}$, maximum value of the Laplacian of the cutoff function used in the proof of \cref{equation:schauder-estimate-on-T3}
\item $\Cr{D-chi-bound}$, maximum value of the derivative of the cutoff function. The cutoff function is defined in \cref{proposition:cut-off-estimates}.
\item $\Cr{D^2-chi-bound}$, largest entry of the second derivative of the cutoff function, \cref{proposition:cut-off-estimates}
\item $\Cr{Holder Q-Qtilde} = 27^\frac{3}{4}$, This is the explicit Holder Constant in the following case. 
\begin{align*}
 \|{f}_{L^4([-1,2]^3)} \leq \Cr{Holder Q-Qtilde} \|{f}_{L^6([-1,2]^3)}
\end{align*}
\end{itemize}
\end{subsection}
\begin{subsection}{\cref{subsection:estimates-on-t3-with-non-flat-metric}}
\begin{itemize}
\item $\Cr{const:metric-max-difference}$: denotes an amount we are able to perturb the metric so that our injectivity estimate still holds, see \cref{theorem:T3-injectivity-estimate}
\item $\Cr{const:injectivity-main-nonflat}$, The main result of Section 3.2, see \cref{theorem:T3-injectivity-estimate}
\item $\Cr{Step 1 3.3}$ The explicit upper bound in \cref{Step 1: T3 Injectivity Estimate}
\item $\Cr{Step 3 3.3}$ \cref{overall strategy for injective estimate with arbitrary-ish riemannian metric} 
\item $\Cr{Step 1 3.3 Lower Bound}$ See \cref{Step 1: T3 Injectivity Estimate}, this is an estimate that is required for the application in Section 4.
\item $\Cr{Lower Determinant Estimate} = (1-\delta)^3- 2\delta^3 -3(1+\delta) \delta^2$ A lower bound for the determinant of the perturbed metric, see \cref{Determinant Estimate}
\item $\Cr{Upper Determinant Estimate} = (1+\delta)^3+ 2\delta^3 +3(1+\delta) \delta^2$ An upper bound for the determinant of the perturbed metric, see \cref{Determinant Estimate}
\item $\Cr{Covector Lower Estimate}$ Compares pointwise norms of covectors when taken with respect to the flat metric and the perturbed metric, this is the lower bound, see \cref{comparison lemma covectors}
\item $\Cr{Covector Estimate}$ Compares pointwise norms of covectors when taken with respect to the flat metric and the perturbed metric, this is the upper bound, see \cref{comparison lemma covectors}
\item $\Cr{Co2Tensor Estimate}$ Compares pointwise norms of covariant 2-tensors when taken with respect to the flat metric and the perturbed metric, see \cref{comparison lemma covectors}
\item $\Cr{Christoffel Bound}$ See \cref{lemma:christoffel-bound}, this is a bound on the value of the Christoffel Symbols for the perturbed metric.
\item $\Cr{S1P1 T3 Injectivity Estimate Lower Bound}$ is an explicit value for the lower bound in \cref{Step 1 Part 1: T3 Injectivity Estimate}
\item $\Cr{S1P1 T3 Injectivity Estimate}$ is an explicit value for the upper bound in \cref{Step 1 Part 1: T3 Injectivity Estimate}, 
\item $\Cr{Absorption Laplacian Term}$, the constant computed in \cref{Delta Comparison Bound}, this is necessary to compare the norm of the laplacian between the flat and perturbed metrics. 
\end{itemize}
\end{subsection}
\begin{subsection}{\cref{section:application-nowhere-vanishing-1-forms}}
\begin{itemize}
\item $\Cr{endofpaper}$ is defined in \cref{End of Proof}
\item $\Cr{sobolev-embedding}$ The optimal Sobolev Embedding constant for the embedding $L^2_{2,0} \hookrightarrow L^4$, this is \cref{Sobolev Embedding}
\item $\Cr{Sobolev-Embedding-Cube}$, an explicit Sobolev Embedding constant for the embedding $L^2_2 \hookrightarrow{L^4} $ when $\Omega=[0,1]^3$, this is \cref{other-sobolev-embedding}
\item $\Cr{Morrey}$ An explicit bound for the constant appearing in Morrey's Inequality. 
\end{itemize}
\end{subsection}
\bibliographystyle{apalike}
\bibliography{library}
\end{document}

%% file: 1forms_preamble.tex
\usepackage[margin=3.2cm]{geometry}

\usepackage{setspace}

\usepackage{graphics}

\usepackage{docmute}

\usepackage[utf8]{inputenc}

\usepackage{mathtools}

\usepackage{amsmath}

\usepackage{amsfonts}
\usepackage{amssymb}

\usepackage{amsthm}
\usepackage{aliascnt}

\usepackage{makeidx}
\makeindex

\usepackage{textcomp}
\usepackage{lmodern}
\usepackage[varqu,varl]{zi4}%
\usepackage{newpxmath}
\usepackage{superiors}
\usepackage{bm}

\usepackage[breaklinks,linktocpage]{hyperref}
\usepackage[hyperpageref]{backref}

\usepackage{url}
\usepackage{breakurl}

\usepackage{tikz-cd}

\usepackage{enumitem}
\usepackage[UKenglish]{babel}
\usepackage[makeroom]{cancel}
\usepackage{pstricks-add}

\definecolor{shadecolor}{rgb}{0.88,0.91,0.95}       
\usepackage{tcolorbox}

\usepackage{dynkin-diagrams}
\usetikzlibrary{backgrounds}

\usepackage{booktabs}

\usepackage{fancyhdr}
\newcommand{\R}{\mathbb{R}}

\newcommand{\Ker}{\operatorname{Ker}}

\renewcommand{\|}[1]{\left| \left| #1 \right| \right|}

\newcommand{\gflat}{g_{\mathrm{flat}}}

\newcommand{\vol}{\operatorname{vol}}

\renewcommand{\d}{\operatorname{d} \!}

\renewcommand{\tilde}[1]{\widetilde{#1}}

\newcommand\Item[1][]{%
  \ifx\relax#1\relax  \item \else \item[#1] \fi
  \abovedisplayskip=0pt\abovedisplayshortskip=0pt~\vspace*{-\baselineskip}}

\usepackage[capitalise]{cleveref}

\numberwithin{equation}{section}

\newtheorem{proposition}[equation]{Proposition}
\crefname{proposition}{Proposition}{Propositions}

\newtheorem{lemma}[equation]{Lemma}
\crefname{lemma}{Lemma}{Lemmas}

\newtheorem{corollary}[equation]{Corollary}
\crefname{corollary}{Corollary}{Corollaries}

\newtheorem*{corollary*}{Corollary}
\crefname{corollary*}{Corollary}{Corollaries}

\newtheorem{theorem}[equation]{Theorem}
\crefname{theorem}{Theorem}{Theorems}

\crefname{conjecture}{Conjecture}{Conjectures}

\newtheorem*{theorem*}{Theorem}
\crefname{theorem*}{Theorem}{Theorems}

\crefname{claim}{Claim}{Claims}

\theoremstyle{remark}

\crefname{question}{Question}{Questions}

\newtheorem{definition}[equation]{Definition}
\crefname{definition}{Definition}{Definitions}

\crefname{example}{Example}{Examples}

\newtheorem{remark}[equation]{Remark}
\crefname{remark}{Remark}{Remarks}

\crefname{assumption}{Assumption}{Assumptions}